\RequirePackage{fix-cm}
\documentclass[smallextended]{svjour3}       
\smartqed  

\usepackage{lipsum}
\usepackage{graphicx}
\usepackage{multirow}
\usepackage{subfigure}
\usepackage{tikz}
\usepackage{pgfplots} 
\pgfplotsset{compat=1.15}
\usepackage{lineno,hyperref}
\usepackage{algorithm}
\usepackage{algpseudocode}
\usepackage{amsmath}
\usepackage{amsfonts}
\usepackage{amssymb}
\usepackage{subfigure}
\usepackage{overpic}

\begin{document}

\title{An Alternating Direction Implicit Method for Mean Curvature Flows
}


\author{Han Zhou         \and
        Shuwang Li   \and
        Wenjun Ying
}



\institute{
        H. Zhou 
        \at School of Mathematical Sciences, Shanghai Jiao Tong University, Minhang, Shanghai 200240, PR China \\
        \email{zhouhan@sjtu.edu.cn}
        \and S. Li 
        \at Department of Applied Mathematics, Illinois Institute of Technology, Chicago, USA\\ 
        \email{sli15@iit.edu} 
        \and 
        W. Ying 
        \at Corresponding author. School of Mathematical Sciences, MOE-LSC and Institute of Natural Sciences, Shanghai Jiao Tong University, Minhang, Shanghai 200240, PR China\\ 
        \email{wying@sjtu.edu.cn}}

\date{Received: date / Accepted: date}

\maketitle






\begin{abstract}
This paper is concerned with the mean curvature flow, which describes the dynamics of a hypersurface whose normal velocity is determined by local mean curvature.
We present a Cartesian grid-based method for solving mean curvature flows in two and three space dimensions.
The present method embeds a closed hypersurface into a fixed Cartesian grid and decomposes it into multiple overlapping subsets.
For each subset, extra tangential velocities are introduced such that marker points on the hypersurface only moves along grid lines.
By utilizing an alternating direction implicit (ADI)-type time integration method, the subsets are evolved alternately by solving scalar parabolic partial differential equations on planar domains.
The method removes the stiffness using a semi-implicit scheme and has no high-order stability constraint on time step size.
Numerical examples in two and three space dimensions are presented to validate the proposed method.
\end{abstract}
\section{Introduction}\label{sec:intro}
Geometric evolution of interfaces draws many attentions in the last decades due to its wide applications in mathematics \cite{Colding2015}, materials science \cite{Mullins1956}, biology \cite{Nelson2005} and, more recently, image processing \cite{Malladi1995,6253258,Mallat2007}.
In a geometric evolution problem, the dynamics of a hypersurface are described by its geometry.
Typically, the normal velocity of the hypersurface is given by a law defined by geometry.
In this paper, we are concerned with a representative case of geometric evolution problems, the mean curvature flow, in which the hypersurface evolves such that its normal velocity equals its negative mean curvature.
Mean curvature flow was originally proposed by Mullins to model an ideal grain boundary motion \cite{Mullins1956}.
Thereafter, it was also used to model various other physical phenomena \cite{Nelson2005,Alias2020}.

Numerical methods for mean curvature flows can be classified into three categories, based on their different representations of hypersurfaces, which are parametric approaches, level set method \cite{osher2004level,Osher1988,osher2003geometric} and phase field method \cite{Bellettini1996,Chen1992,paolini1995efficient}.
Representative parametric approaches include the parametric finite element method \cite{Bao2021,Barrett2020,BARRETT20084281,BARRETT2007441}, graph approach \cite{deckelnick_dziuk_1999,Deckelnick2005} and the front tracking method \cite{Lai2009,UNVERDI199225}.
For a comprehensive survey on numerical methods for mean curvature flows, the interested reader is referred to the review article by Deckelnick et al. \cite{Deckelnick2005}.

For general moving interface problems, although the level set method and the  phase field method have their advantages in handling topological changes and ease in implementation, parametric approaches provide surprisingly good results, such as accurate computation of curvature and conservation of mass, even with a coarse grid, and are computationally very cheap as well \cite{Barrett2014,Tryggvason2001}.
Despite the benefits, numerical computation of mean curvature flows with parametric approaches also encounters several difficulties, including the deterioration of mesh quality during the computation and the numerical stiffness induced by the mean curvature term.
Due to the pure normal motion of the surface, adjacent mesh nodes may become closer and closer, making the computation highly unstable.
The problem is even more severe for the mean curvature flow due to its "curve shortening" property by Mullins \cite{Mullins1956}.
In addition, the evolution equation of mean curvature flows has second-order spatial derivatives in the mean curvature term, which induces numerical stiffness such that a small time step is required for explicit time integration schemes \cite{HOU1994312}.
A naive discretization with implicit time integration to remove stiffness leads to a nonlinear system, for which finding a numerical solution is time-consuming, even with advanced iterative solvers.
For two space dimensional curves, these difficulties can be very well handled by the small scale decomposition method, initially proposed by Hou et al. \cite{HOU1994312}.
The idea was also extended to three space dimensional cases for some special surfaces \cite{HOU1998628,Ambrose2013}.
However, for general closed surfaces in three space dimensions, it is still unclear how to apply small scale decomposition for mean curvature flows.

This work proposes a new numerical method for mean curvature flows.
The method decomposes a moving hypersurface into multiple overlapping subsets such that each subset can be viewed as a Monge patch for which the graph approach \cite{deckelnick_dziuk_1999,Deckelnick2005} is applicable.
The method is based on an overlapping decomposition method for hypersurfaces, initially proposed by Wilson for computing integrals on implicitly defined curves and surfaces \cite{Jason2010}.
A few years later, the method was extended by the second author to incorporate a kernel-free boundary integral method for solving elliptic partial differential equations on irregular domains \cite{Ying2013}.
The decomposition strategy has many advantages in simplicity and efficiency.
By representing each subset with its intersection points with grid lines, it is natural to keep marker points quasi-equidistant and maintain mesh quality.
By reformulating the evolution equation, which is a nonlinear system, into a sequence of scalar PDEs 
on overlapping subsets, one can evolve the subsets alternately in the spirit of the alternating direction implicit (ADI) method \cite{Peaceman1955,Douglas1962,Douglas1964}.
The resulting algorithm is efficient and has no high-order stability constraint.

The remainder of the paper is organized as follows. In Section \ref{sec:math}, we describe the governing equation of mean curvature flow and its hybrid formulation based on an overlapping surface decomposition method. The numerical methods for solving mean curvature flows are described in Section \ref{sec:method}. In Section \ref{sec:algorithm}, the numerical algorithm of the proposed method is briefly summarized. Multiple numerical examples are presented to validate the present method in Section \ref{sec:result}. In the final Section \ref{sec:discu}, we briefly discuss the present method and some further work.
\section{Mathematical formulation}\label{sec:math}

Let $\Gamma(t)\subset\mathbb{R}^d, d =2, 3$ be a closed moving hypersurface. 
Consider the mean curvature flow problem that, for any point $\boldsymbol{x}$ on $\Gamma$, the evolution is given by,
\begin{equation}\label{eqn:mcf}
\boldsymbol{x}_t = V \boldsymbol{n}, \quad V = -\kappa, \quad  \boldsymbol x \in \Gamma,
\end{equation}
where $\kappa$ is the (mean) curvature and $\boldsymbol{n}$ the unit outward normal. 
Here, a circle/sphere has positive curvature.
By applying the transport theorem of evolving hypersurfaces, it can be shown that the mean curvature flow has length-decreasing and area-decreasing properties in 2D and 3D, respectively,
\begin{align}
    \dfrac{d}{dt}|\Gamma(t)| &= \dfrac{d}{dt}\int_{\Gamma(t)} 1 \,ds =  - \int_{\Gamma(t)} \kappa^2 \, ds < 0, &\text{ if } d = 2,\\
    \dfrac{d}{dt}|\Gamma(t)| &= \dfrac{d}{dt}\int_{\Gamma(t)} 1 \,dA =  - \int_{\Gamma(t)} \kappa^2 \, dA < 0, &\text{ if } d = 3,
\end{align}
where $|\Gamma(t)|$ is the length of $\Gamma(t)$ in 2D and the area of $\Gamma(t)$ in 3D.
Specially, for the case $d=2$, it holds that
\begin{equation}
    \dfrac{d}{dt} A(t) = \int_{\Gamma(t)} \dfrac{d\boldsymbol{x}}{dt} \cdot \boldsymbol{n}\, ds = -\int_{\Gamma(t)} \kappa\,ds = -2\pi,
\end{equation}
where $A(t)$ is the area enclosed by $\Gamma(t)$.
It suggests that the enclosed area of a 2D mean curvature flow always decreases at a constant rate \cite{Dallaston2016}.
The result no longer holds for the enclosed volume of 3D surfaces since the surface integral of mean curvature varies from case to case.

If $\Gamma$ is closed, and there is no boundary condition, then the solution of mean curvature flow is determined by its initial configuration.
For some certain initial configurations, solutions of mean curvature flows may develop singularities, such as pinch-off and topological changes, in finite time before shrinking to a point.
In this paper, we seek well-defined solutions of mean curvature flows, $i.e.$, solutions before singularities happen.

In order to tackle the mean curvature flow problem, the evolution equation \eqref{eqn:mcf} is divided into multiple subproblems with an overlapping surface decomposition of $\Gamma$. 
For $r = 1,\cdots, d$, let $\mathbf{e}_r$ be the $r^{th}$ unit vector in $\mathbb{R}^d$ and $\alpha\in(\cos^{-1}(1/\sqrt{d}),  \pi/2)$ be a fixed angle. The set 
\begin{equation}\label{eqn:sur-decom}
\Gamma_r = \{\boldsymbol{x}\in \Gamma :| \mathbf{n}\cdot\mathbf{e}_r|(\boldsymbol{x})>\cos\alpha\},
\end{equation}
is an open subset of the surface $\Gamma$ for each $r = 1, 2, \cdots, d$. The union of the sets $\{\Gamma_r\}_{r = 1}^{d}$ forms an overlapping surface decomposition of $\Gamma$. 
Then the surface $\Gamma$ is represented by the overlapping subsets $\Gamma_r$ with a partition of unity.

\subsection{Divided problems}
Note that the evolution of a  hypersurface is only determined by its normal velocity $V$. 
One is allowed to add arbitrary tangential velocity to the evolution of $\Gamma$ without altering its shape. The tangent velocity only changes in the frame for the parametrization of the surface. 
For arbitrary tangential velocities $T$, $T_1$ and $T_2$, the evolution governed by \eqref{eqn:mcf} is equivalent to
\begin{equation}\label{eqn:mcf-2d-mod}
\boldsymbol{x}_t = -V\pmb{n} + T \pmb{\tau},\quad \boldsymbol{x}\in\Gamma,
\end{equation}
in two space dimensions and
\begin{equation}\label{eqn:mcf-3d-mod}
\boldsymbol{x}_t = -V\pmb{n} + T_1 \pmb{\tau}_1 + T_2 \pmb{\tau}_2,\quad \boldsymbol{x}\in\Gamma,
\end{equation}
in three space dimensions.
Here, the notations $\boldsymbol{\tau}$, $\boldsymbol{\tau}_1$ and $\boldsymbol{\tau}_2$ mean tangent vectors.

Consider the evolution of the overlapping subsets $\Gamma_r$ in three space dimensions.
Let $\boldsymbol x^{(r)}$ denote a point on the subset $\Gamma_r$. 
The evolutions of $\Gamma_r$ are given by 
\begin{equation}
\boldsymbol{x}^{(r)}_t = - V^{(r)}\boldsymbol{n}, \quad \boldsymbol{x}^{(r)}\in\Gamma_r , \quad r = 1, \cdots, d.
\end{equation}
where $V^{(r)}$ are the restrictions of $V$ from $\Gamma$ to $\Gamma_r$.
By adding tangential velocities $T^{(r)}_1$ and $T^{(r)}_2$, it yields the equivalent evolution equations
\begin{equation}\label{eqn:decom-form}
\boldsymbol{x}^{(r)}_t = - V^{(r)} \boldsymbol{n} + T^{(r)}_1 \boldsymbol{\tau}_1 + T^{(r)}_2 \boldsymbol{\tau}_2, \quad \boldsymbol{x}^{(r)}\in\Gamma_r , \quad r = 1, \cdots, d.
\end{equation}
Hence, the evolution equation \eqref{eqn:mcf} of $\Gamma$ is divided into a sequence of evolution equations of subsets $\Gamma_r$.
With the divided formulation \eqref{eqn:decom-form}, for each subset $\Gamma$, the tangential velocities $T^{(r)}_1$ and $T^{(r)}_2$ can be chosen independently.

Due to the overlapping surface decomposition \eqref{eqn:sur-decom}, each subset $\Gamma_r$ can be easily parameterized with Cartesian coordinates in a planar domain $\Omega_r\subset \mathbb{R}^{d-1}$. 
For example, in three space dimensions, denote by $\Omega_3$ the projection of $\Gamma_3$ onto $X$-$Y$ plane. 
Then $\Gamma_3$ can be represented by the Monge patch $\boldsymbol{x}(x,y) = \boldsymbol{x}(x,y, z(x,y)),(x,y)\in\Omega_3$ in which $z(x, y)$ is a height function. 
With this understanding, the evolution of $\Gamma_r$ can be described as a time-dependent height function on the base plane $\Omega_3$ in $d-1$ space dimensions. 
The height function representation is an Eulerian description of the moving hypersurface.
Its numerical approximation is much simpler than that for tracking a moving hypersurface with its Lagrangian motion.
With Eulerian description, it is natural to use a Cartesian grid to approximate the height function.
With the understanding that the Eulerian description is equivalent to moving marker points of the hypersurface along fixed grid lines, the evolution equations of the equivalent Eulerian motion can be derived by carefully choosing tangent velocities $T^{(r)}$, $T^{(r)}_1$ and $T^{(r)}_2$ such that $\Gamma_r$ only have one non-zero velocity component in direction $\boldsymbol{e}_r$.

\subsubsection{Two space dimensional case}
Let $\Gamma\subset \mathbb{R}^2$ be a time-dependent Jordan curve which is given by $\boldsymbol{x}(t) = (x(\theta, t), y(\theta, t))$ where $\theta$ parameterizes the curve. Its curvature $\kappa$ and unit outward normal vector are, respectively, given by
\begin{equation}\label{eqn:cur-normal}
\kappa = \frac{x_{\theta}y_{\theta\theta} - x_{\theta\theta}y_{\theta}}{(x_{\theta}^{2} + y_{\theta}^{2})^{\frac{3}{2}}},\quad \boldsymbol{n} = \dfrac{1}{(x_{\theta}^{2} + y_{\theta}^{2})^{\frac{1}{2}}} \begin{pmatrix}
y_{\theta}\\
-x_{\theta}
\end{pmatrix}.
\end{equation}

For each subset $\Gamma_r, r = 1,2$, it can be parameterized by $x$ or $y$ to be a height function $y = y(x)$ or $x = x(y)$ depending on its orientation. 
Suppose $\Gamma_r$ is represented in the form $\eta = \eta(\xi)$ where $(\xi, \eta)$ coincides with $(x,y)$ or $(y, x)$.
After extra tangential velocity $T$ is added into the original evolution equation \eqref{eqn:mcf}, the evolution of $\Gamma_r$ is equivalent to
\begin{align}\label{eqn:par-form-2d}
\frac{d}{dt}
\begin{pmatrix}
\xi\\
\eta
\end{pmatrix}
= &- \frac{\eta_{\xi\xi}}{(\eta_{\xi}^{2} + 1)^{2}}
\begin{pmatrix}
\eta_{\xi}\\
- 1
\end{pmatrix}
 + T
 \begin{pmatrix}
1\\
\eta_{\xi}
\end{pmatrix}.
\end{align}
{\color{black}It is worth mentioning that equation \eqref{eqn:par-form-2d} does not rely on the orientation of the curve due to the cancellation of signs in the curvature and normal vector when one reverses the parameterization, namely, from $\xi$ to $-\xi$.}

To determine a specific tangential velocity $T$ such that marker points on $\Gamma_r$ only have non-zero velocity component in $\boldsymbol{e}_r$ direction.
One needs to set $\xi_t = 0$, $i.e.$
\begin{equation}
\xi_t = -\frac{\eta_{\xi}\eta_{\xi\xi}}{(\eta_{\xi}^{2} + 1)^{2}} + T = 0.
\end{equation}
The expression of $T$ can be easily solved.
By substituting the determined $T$ into \eqref{eqn:par-form-2d}, it yields the evolution law for $\Gamma_r$ in terms of height function, 
\begin{equation}\label{eqn:2d-pde}
\eta_t = \frac{\eta_{\xi\xi}}{\eta_{\xi}^{2} + 1},
\end{equation}
which is a scalar parabolic-type partial differential equation.

\subsubsection{Three space dimensional case}
The derivation in three space dimensions is similar. 
For each subset $\Gamma_r, r = 1,2,3$, it can be regarded as a Monge patch $\boldsymbol{x}(u,v) = \boldsymbol{x}(u,v,w(u, v)), (u, v)\in\Omega_r$ where $\Omega_r$ is the projection of $\Gamma_r$ on its base plane and $w(u,v)$ is the height function.
Denote by $\boldsymbol{\tau}_1 = (1, 0, w_u)^T,\boldsymbol{\tau}_2 = (0, 1, w_v)^T$ two tangent vectors of $\Gamma_{r}$. 
After adding two tangential velocities $T_1$ and $T_2$, the evolution equation $\eqref{eqn:mcf}$ in three space dimensions becomes
\begin{equation}\label{eqn:par-3d}
\frac{d}{dt}
\begin{pmatrix}
u\\
v\\
w
\end{pmatrix} =  
\frac{(1 + w_u^2)w_{vv} - 2w_{u}w_{v}w_{uv} + (1 + w_{v}^2)w_{uu}}{2(1 + w_{u}^2 + w_{v}^2)^2}
\begin{pmatrix}
-w_{u}\\
-w_v\\
1
\end{pmatrix} + 
T_1
\begin{pmatrix}
1\\
0\\
w_u
\end{pmatrix} + 
T_2
\begin{pmatrix}
0\\
1\\
w_v
\end{pmatrix}.
\end{equation}
By setting $u_t = v_t = 0$, one can solve for $T_1$ and $T_2$,
\begin{align}
T_1 = w_u \frac{(1 + w_u^2)w_{vv} - 2w_{u}w_{v}w_{uv} + (1 + w_{v}^2)w_{uu}}{2(1 + w_{u}^2 + w_{v}^2)^2}\label{eqn:tan-v1}, \\
T_2 = w_v \frac{(1 + w_u^2)w_{vv} - 2w_{u}w_{v}w_{uv} + (1 + w_{v}^2)w_{uu}}{2(1 + w_{u}^2 + w_{v}^2)^2}\label{eqn:tan-v2}.
\end{align}
By substituting \eqref{eqn:tan-v1} and \eqref{eqn:tan-v2} into \eqref{eqn:par-3d}, it gives the evolution of $\Gamma_{r}$ in terms of its height functions $w$,
\begin{equation}\label{eqn:3d-pde}
w_t = \frac{(1 + w_u^2)w_{vv} - 2w_{u}w_{v}w_{uv} + (1 + w_{v}^2)w_{uu}}{2(1 + w_{u}^2 + w_{v}^2)}.
\end{equation}
The equation \eqref{eqn:3d-pde} is also a scalar parabolic-type partial differential equation.

\subsection{Matching condition}
Until now, the original evolution equation $\eqref{eqn:mcf}$ is reformulated as a sequence of scalar partial differential equations $\eqref{eqn:2d-pde}$ and $\eqref{eqn:3d-pde}$. 
To ensure the well-posedness of the divided problem, we follow the idea of domain decomposition \cite{mathew2008domain} to add an extra matching condition for the solution of the divided problem at the overlapping zone such that the equations \eqref{eqn:2d-pde} and \eqref{eqn:3d-pde} have boundary condition.

A simple choice of the matching condition is to enforce continuity of the global solution at an overlapping zone with a partition of unity,
\begin{equation}\label{eqn:match-cond}
\boldsymbol{x}^{(r)} = \sum_{j \neq r}\chi_j \boldsymbol{x}^{(j)}, \quad \boldsymbol{x}^{(r)}\in \partial\Gamma_r.
\end{equation}
where $\partial\Gamma_r$ denotes the boundary of $\Gamma_r$. 
Here, the notation $\chi_j$ is the partition of unity subordinate to the subset $\Gamma_j$, which satisfies
\begin{equation}
\left \{
\begin{aligned}
&\chi_r(\boldsymbol{x}) \geq 0, \quad &\boldsymbol{x} \in \Gamma_r,\\
&\chi_r(\boldsymbol{x}) = 0, \quad &\boldsymbol{x} \in \Gamma \backslash \Gamma_r, \\ 
&\sum_{r=1}^d\chi_r(\boldsymbol{x}) = 1, \quad &\boldsymbol{x} \in \Gamma.
\end{aligned}
\right.
\end{equation}

The divided problem \eqref{eqn:decom-form} together with the matching condition \eqref{eqn:match-cond} forms an equivalent coupled system, which is called hybrid formulation, to equation \eqref{eqn:mcf} in three space dimensions,
\begin{equation}\label{eqn:hybrid}
\left \{
\begin{aligned}
&\boldsymbol{x}^{(r)}_t = - V^{(r)} \boldsymbol{n} + T^{(r)}_1 \boldsymbol{\tau}_1 + T^{(r)}_2 \boldsymbol{\tau}_2, \quad &\boldsymbol{x}^{(r)}\in\Gamma_r, \\
&\boldsymbol{x}^{(r)} = \sum_{j \neq r}\chi_j \boldsymbol{x}^{(j)}, \quad &\boldsymbol{x}^{(j)}\in \partial\Gamma_r.
\end{aligned}
\right.
\end{equation}
Once the hybrid formulation \eqref{eqn:hybrid} is solved, the global solution $\boldsymbol{x}$ can be reconstructed with a partition of unity, 
\begin{equation}
\boldsymbol{x} = \sum_{r = 1}^{d}\chi_r\boldsymbol{x}^{(r)}.
\end{equation} 
Unlike Ambrose's method \cite{Ambrose2013}, which is applicable only for a particular class of surfaces with doubly-periodic boundary conditions, our method can handle more general cases, including closed surfaces,  with this hybrid formulation.
\section{Numerical Methods}\label{sec:method}
In this section, the numerical methods for mean curvature flow are described, including the discrete representation of a moving hypersurface, numerical discretizations of the partial differential equations $\eqref{eqn:2d-pde}$ and $\eqref{eqn:3d-pde}$ as well as the matching condition $\eqref{eqn:match-cond}$.

\subsection{Hypersurface representation}\label{sec:osd-dis}
\begin{figure}[htbp]
\centering
\subfigure[]{\includegraphics[width=0.4\textwidth]{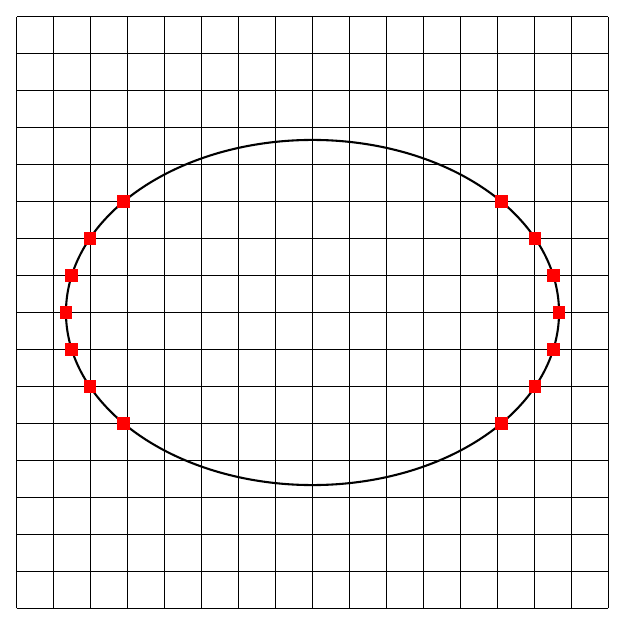}}
\subfigure[]{\includegraphics[width=0.4\textwidth]{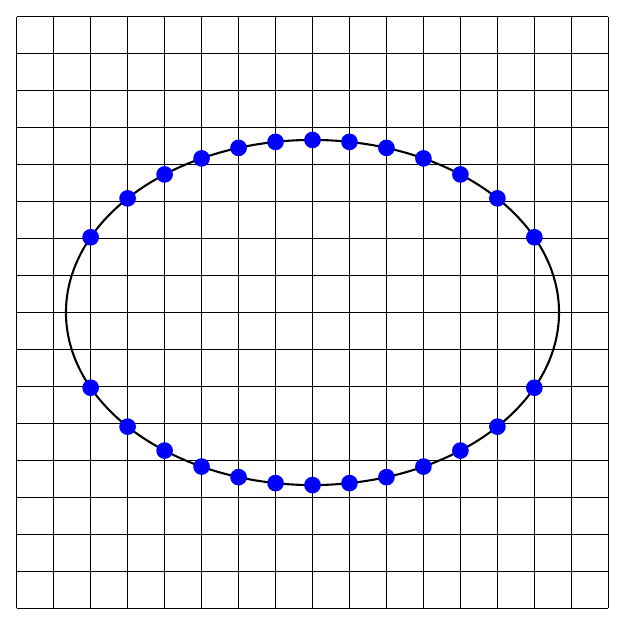}}
\caption{Control points for the representation of an ellipse: (a) control points on $\Gamma_1$ (red rectangle markers); (b) control points on $\Gamma_2$ (blue circle markers).}
\label{fig:osd-nodes}
\end{figure}

\begin{figure}[htbp]
\centering
\includegraphics[width=0.5\textwidth]{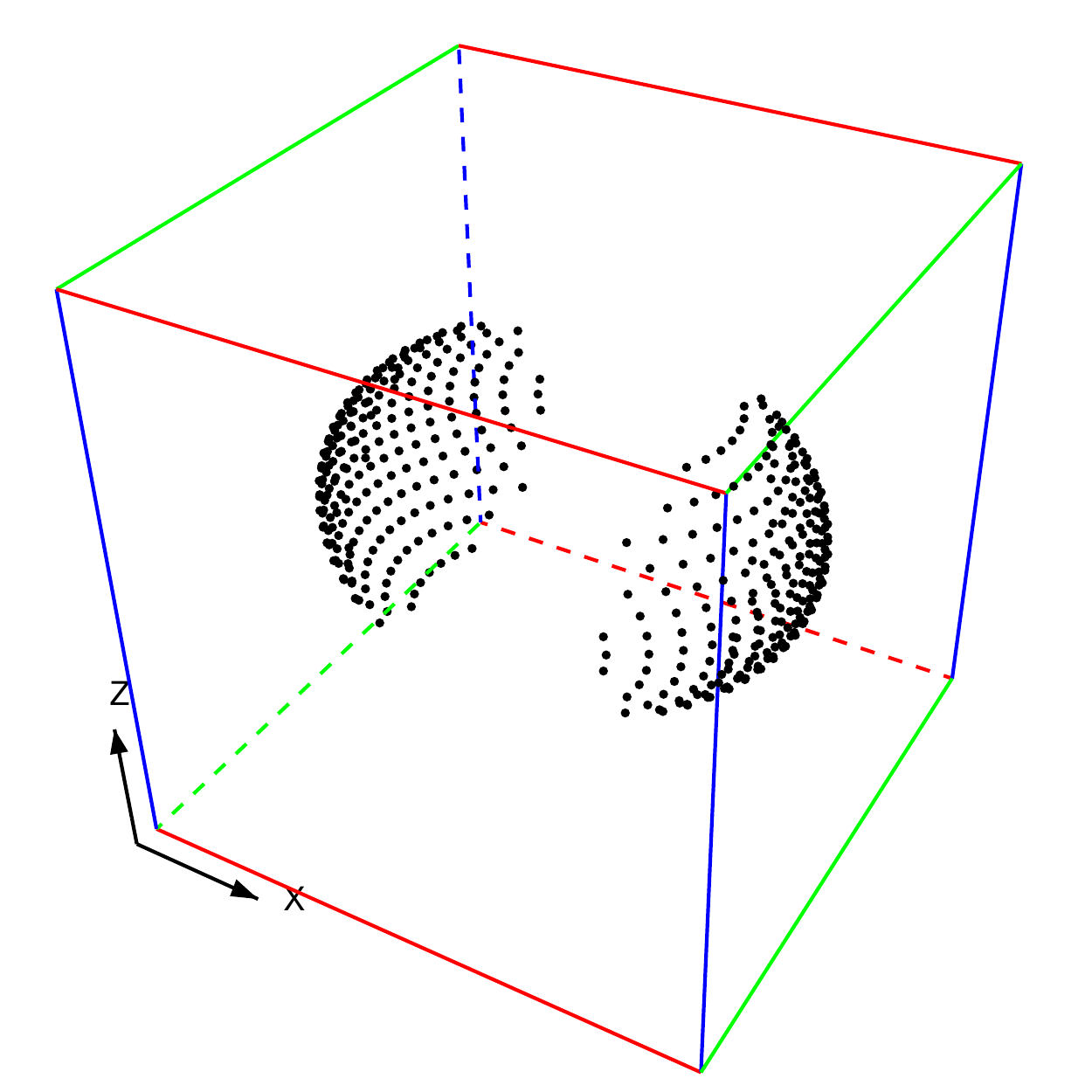}
\caption{Control points for the representation of an ellipsoid on $\Gamma_1$.}
\label{fig:osd-nodes-3d}
\end{figure}

Let $\Gamma$ be a smooth closed hypersurface. 
It is separately represented by the height functions of its overlapping subsets $\Gamma_r, r = 1,\cdots,d$, which are approximated by nodal values at Cartesian grid nodes in the base domain $\Omega_r$.
Equivalently, the nodal values are, in fact, the intersection points of $\Gamma_r$ and grid lines, which are aligned with $\boldsymbol{e}_r$.
At the implementation level, this is done by selecting from all intersection points  $\boldsymbol{p}$ of $\Gamma$ and grid lines for certain ones which satisfy the decomposition rule:
\begin{equation}\label{eqn:decom-rule}
     \boldsymbol p\in\Gamma_r\text{ and } \boldsymbol{n}(\boldsymbol{p})\cdot \boldsymbol{e}_r>\cos\alpha.
\end{equation}
where $\boldsymbol{n}(\boldsymbol{p})$ is the unit outward normal at $\boldsymbol p$ and $\alpha$ is a given threshold angle.
Those points which satisfy \eqref{eqn:decom-rule} are named as control points for $\Gamma_r$ and the point set is denoted by $\Gamma^h_r$.
We also denote all control points on $\Gamma$ by $\Gamma^h = \cup_{r=1}^d \Gamma^h_r$.
The point set $\Gamma^h$ is used to represent $\Gamma$ in terms of its overlapping subsets.
We remark that, although points in $\Gamma^h$ are not quasi-equidistant, local interpolation stencil only involves points in each subset $\Gamma^h_r$, which are quasi-equidistant.
Figure \ref{fig:osd-nodes} and \ref{fig:osd-nodes-3d} show the distribution of control points in two and three space dimensions, respectively.


The advantage of this representation of hypersurface is evident.  
One can easily find out that the projections of control points in $\Gamma^h_r$ coincide with Cartesian grid nodes in $\Omega_r$.
Instead of tracking $\Gamma$ by marker points whose velocities are in $d$ dimensions, one needs to solve the evolution for the height functions, $i.e.$ marker points moving along grid lines, which only change values in one dimension.



\subsubsection{Solving PDEs on hypersurfaces}

Generally, for a closed smooth hypersurface $\Gamma$, the subset $\Gamma_r$ consists of several isolated components which are denoted by $\Gamma_{r,l}, l=1,2,\cdots$. 
For example, suppose $\Gamma:\boldsymbol{x} = \boldsymbol{x}(\theta),\theta\in[0,2\pi)$ is a circle, then the curve segments $\Gamma_{1,1}:\boldsymbol{x} = \boldsymbol{x}(\theta),\theta\in(2\pi-\alpha,2\pi)\cup[0, \alpha)$ and $\Gamma_{1,2}:\boldsymbol{x} = \boldsymbol{x}(\theta),\theta\in(\pi-\alpha, \pi+\alpha)$ are both subsets of $\Gamma_1$. 
Let $\Omega_{r,l}$ denote the projection of $\Gamma_{r,l}$ on the base plane.
If $\Omega_{r,l}$ overlap with each other, then $\Gamma_r$ is a multi-valued function on $\Omega_r$, which induces ambiguity. 

\begin{figure}[htbp]
    \centering
    \includegraphics[width=0.5\textwidth]{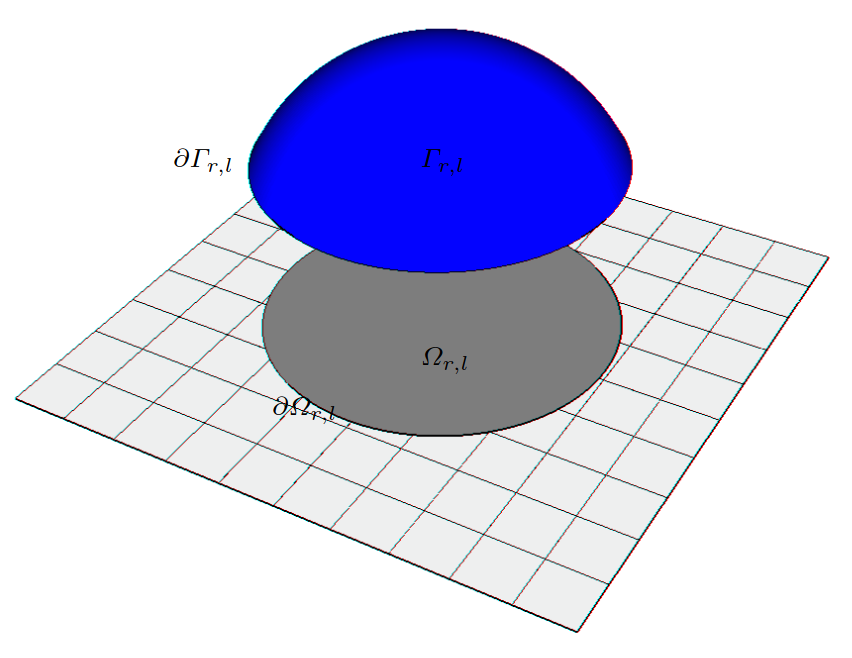}
    \caption{A surface component $\Gamma_{r,l}$ and its projection on the base plane $\Omega_{r,l}$.}
    \label{fig:patch_prj}
\end{figure}

The correct understanding is to separate $\Gamma_{r,l}$ from each other, and on which PDEs are solved independently.
Once the components are separated from each other, the ambiguity is removed, since $\Gamma_{r,l}$ is a single-valued function on $\Omega_{r,l}$ (see Figure \ref{fig:patch_prj}).
Only points on the same component are involved in a local stencil for solving PDEs.
In the implementation, one can check their distance and normals to determine if two points are on the same component $\Gamma_{r,l}$.
We identify two points $\boldsymbol{p}$ and $\boldsymbol{q}$ on the same component if they satisfy
\begin{equation}\label{eqn:criteria}
    \Vert \boldsymbol{p} - \boldsymbol{q}\Vert < D_0\text{ and }\boldsymbol{n}(\boldsymbol{p})\cdot \boldsymbol{n}(\boldsymbol{q}) < \cos(\theta_0),
\end{equation}
where $\Vert\cdot\Vert$ is the Euclidean distance, $D_0$ and $\theta_0$ are threshold values given in advance.
In this work, we set $D_0=5h$ and $\theta_0 = \pi/6$. 
The separation procedure is meant to find the correct finite difference stencil point from possible intersection points that lie on the same grid line.
The implementation based on the criteria \eqref{eqn:criteria} can naturally handle cases with multiple (more than 2) components in each $\Gamma_r$, for example, oscillating curves or surfaces, as long as the grid is fine enough to resolve them.

\subsubsection{Interpolation on hypersurface}\label{sec:interp}
In the discrete representation $\Gamma^h$, geometric quantities and functions on the hypersurface are evaluated by local interpolation.
Given a point $\boldsymbol{p}\in\Gamma$, to interpolate the function value at $\boldsymbol{p}$ using function values on $\Gamma_r^h$, one needs to find its projection point $\boldsymbol{p}^{\star}\in\Omega$ in which the selection of $r$ depends on the direction of $\boldsymbol{n}(\boldsymbol{p})$.
A quadratic polynomial is locally constructed for interpolation,
\begin{equation}
P_2(u,v) = c_1 + c_2 u + c_3 v + c_4 u^2 + c_5 v^2 + c_6 uv,
\end{equation}
where $(u, v)$ is the local coordinate near $\boldsymbol{p}^{\star}$ in the base plane.
With the help of a Cartesian grid, finding interpolation stencils on $\Gamma^h$ is very simple. One can attach control points to their closest grid nodes and find stencil points by searching nearby grid nodes with their indices in the Cartesian grid. 

\subsubsection{Evolving the hypersurface}
After the PDEs \eqref{eqn:2d-pde} and \eqref{eqn:3d-pde} are solved for a time step, control points are moved to different positions and form a new hypersurface. 
Old points can no longer represent the new hypersurface since they do not satisfy the decomposition rule mentioned before.
To represent the new hypersurface, new control points must be added, and some old ones must be deleted.
This is achieved by finding out all intersection points by local interpolation and select for new control points. 
Even in three space dimensions, intersection points can be computed using one-dimensional polynomial interpolations since both stencil points and new intersection points are on the same plane.

{\color{black}Take the surface component $\Gamma_{3, l}$ as an example.
The component is discretized by intersection points with coordinates $(x_i, y_j, \eta_{i,j})$ where $(x_i,y_j)\in\Omega^h_{3,l}$.
In order to find new intersection points of the component with the grid line $\{(x,y,z)|y=y_j, z=z_k\}$, one first needs to identify the intersection interval $(x_{i_0}, x_{i_1})$ by check the side of grid nodes, and then choose three  intersection points $(x_{i_0-1},y_j,\eta_{i_0-1,j})$, $(x_{i_0},y_j,\eta_{i_0,j})$, and $(x_{i_0+1},y_j,\eta_{i_0+1,j})$ to locally construct a 1D quadratic polynomial $z=P_2(x)$.
By solving the equation $P_2(x^*) = z_k$ with either the Newton method or the bisection method, one can obtain the new intersection point $(x^*, y_j, z_k)$.}

After an intersection point is found, following \eqref{eqn:decom-rule}, one can determine whether it is kept or deleted by checking the normal vector, which is also evaluated by locally constructing a parabola.
It is worth mentioning that if a new component $\Gamma_{r,l}$ is too small and does not have enough control points for interpolation, the whole component should be deleted.

\subsection{Discretization of PDEs}
\subsubsection{Temporal discretization}
The equation \eqref{eqn:2d-pde} and \eqref{eqn:3d-pde} are parabolic PDEs and have second-order spatial derivatives. 
Generally, explicit time integration methods, such as the forward Euler method, for PDEs involving high-order spatial derivatives suffer from high-order stability constraints, and one has to use small time steps to solve the equations.
These problems are known as stiff problems.
The stiffness for mean curvature flow is rather severe. 
Take two space dimensional cases as an example. 
If one discretizes $\Gamma$ by uniformly partitioning the parameter $\theta$. 
An explicit time integration suffers from a second order stability constraint $\Delta t\leq C(\min_{\theta} s_{\theta}h)^2$ where $C$ is a constant,  $s$ is arclength, and $h$ is the grid spacing in $\theta$. 
Since 2-D mean curvature flow is also known as the "curve shortening flow", $s_{\theta}h$ decreases with time and results in an even worse situation.
Although implicit methods are unconditionally stable for stiff problems, there is another difficulty in implicit time integration for \eqref{eqn:2d-pde} and \eqref{eqn:3d-pde} since it results in a nonlinear system in each time step, for which finding a solution is highly inefficient.

Notice that equations \eqref{eqn:2d-pde} and \eqref{eqn:3d-pde} are quasi-linear equations.
The source of stiffness comes from the highest-order terms, which are linear in equations \eqref{eqn:2d-pde} and \eqref{eqn:3d-pde}.
The stiffness can be removed by only treating the highest order terms implicitly with lower order terms treated explicitly.
The resulting time integration scheme is semi-implicit in time.
It only requires solving linear systems, which are much more acceptable than nonlinear systems.

Suppose the time interval $[0, T]$ is uniformly partitioned into $0 = t^0<t^1 <\cdots<t^n<\cdots<t^N = T$ with $t^{n+1}-t^n = \Delta t$.
The semi-implicit schemes for equations $\eqref{eqn:2d-pde}$ and $\eqref{eqn:3d-pde}$, respectively, are give by 
\begin{equation}\label{eqn:2d-semi-dis}
\dfrac{\eta^{n+1} - \eta^n}{\Delta t} = \frac{\eta_{\xi\xi}^{n+1}}{(\eta_{\xi}^n)^{2} + 1},
\end{equation}
and
\begin{equation}\label{eqn:3d-semi-dis}
\dfrac{w^{n+1} - w^n}{\Delta t} = \frac{(1 + (w_u^n)^2)w_{vv}^{n+1} - 2w_{u}^n w_{v}^n w_{uv}^{n+1} + (1 + (w_{v}^n)^2)w_{uu}^{n+1}}{2(1 + (w_{u}^n)^2 + (w_{v}^n)^2)}.
\end{equation}
The semi-implicit schemes \eqref{eqn:2d-semi-dis} and \eqref{eqn:3d-semi-dis} are only in semi-discrete forms. 

\subsubsection{Spatial discretization}
Generally, physical properties should be encoded into the discretization of spatial derivatives. 
In the case of mean curvature flow, the evolution of a hypersurface is driven by surface tension, which is diffusive in nature.
Since diffusion comes from all directions, it is preferred to adopt central differences for approximating the spatial derivatives in the mean curvature term.
Though tangential velocities terms may introduce the convection effect, which commonly should be discretized with methods for hyperbolic PDEs, we remark that the convection effect is expected to be small compared to the diffusion effect.
Hence, for simplicity, we approximated all the spatial derivatives in $\eqref{eqn:2d-pde}$ and $\eqref{eqn:3d-pde}$ with central differences.

Suppose $\Gamma\in\mathbb{R}^d, d =2,3$ is embedded into a bounding box $\mathcal{B}$ which is the tensor product of one dimensional intervals $\mathcal{I}_i, i = 1, 2, \cdots,d$. 
If $\mathcal{I}_i$ is uniformly partitioned into a Cartesian grid $\mathcal{G}_i$ for each $i$, then the tensor product of $\mathcal{G}_i$ forms a natural uniform partition of $\mathcal{B}$, which is also a Cartesian grid and is denoted by $\mathcal{G}$. 
Without loss of generality, we assume equal bounding intervals $\mathcal{I}_i = [a,b], i = 1, 2,\cdots,d$, and they are uniformly partitioned into $N$ intervals. Denote by $h = (b-a)/N$ the mesh parameter.

Let $\Delta_{\xi}, \delta_{\xi}^2$ denote central difference quotients,
\begin{equation*}
    \Delta_{\xi} u_i = \dfrac{u_{i+1} - u_{i-1}}{2h} ,\quad
    \delta_{\xi}^2 u_i = \dfrac{u_{i+1}-2u_i+u_{i-1}}{h^2}.
\end{equation*}
The fully discrete form of equation $\eqref{eqn:2d-pde}$ is given by
\begin{equation}\label{eqn:2d-fully-dis}
\frac{u^{n + 1}_i - u^n_i}{\Delta t} = \frac{\delta^2_{\xi} u_i^{n + 1}}{(\Delta_{\xi} u_i^n)^2 + 1}.
\end{equation}

Similarly, introduce the central difference quotients,
\begin{align*}
\Delta_{u} w_{ij} &= \dfrac{w_{i + 1, j} - w_{i - 1, j}}{2h}, \quad
\Delta_{v} w_{ij} = \dfrac{w_{i, j + 1} - w_{i, j - 1}}{2h}, \\
\delta_{uu}^2 w_{ij} &= \dfrac{w_{i + 1, j} + w_{i - 1, j} - 2w_{ij}}{h^2},\quad
\delta_{vv}^2 w_{ij} = \dfrac{w_{i, j + 1} + w_{i, j - 1} - 2w_{ij}}{h^2}.
\end{align*}
The fully discrete form of equation $\eqref{eqn:3d-pde}$ is given by
\begin{equation}\label{eqn:3d-fully-dis}
\frac{w^{n + 1}_{ij} - w^n_{ij}}{\Delta t}  = C_{ij, 1}^{n}\delta_{vv}^2w_{ij}^{n + 1} + C_{ij, 2}^{n}\Delta_u \Delta_v w_{ij}^{n + 1} + C_{ij, 3}^{n}\delta_{uu}^2w_{ij}^{n + 1},
\end{equation}
where the coefficients $C_{ij,1}$, $C_{ij,2}$ and $C_{ij,3}$ are specified by
\begin{equation*}
\begin{aligned}
C_{ij, 1}^{n} &= \dfrac{1 + (\Delta_{u}w_{ij}^n)^2}{2(1 + (\Delta_{u}w_{ij}^n)^2 + (\Delta_{u}w_{ij}^n)^2)},\\
C_{ij, 2}^{n} &= \dfrac{- \Delta_{u}w_{ij}^n\Delta_{v}w_{ij}^n}{1 + (\Delta_{u}w_{ij}^n)^2 + (\Delta_{u}w_{ij}^n)^2}, \\
C_{ij, 3}^{n} &= \dfrac{1 + (\Delta_{v}w_{ij}^n)^2 }{ 2(1 + (\Delta_{u}w_{ij}^n)^2 + (\Delta_{u}w_{ij}^n)^2)}.
\end{aligned}
\end{equation*}

\begin{figure}[htbp]
\centering
\includegraphics[width=0.4\textwidth]{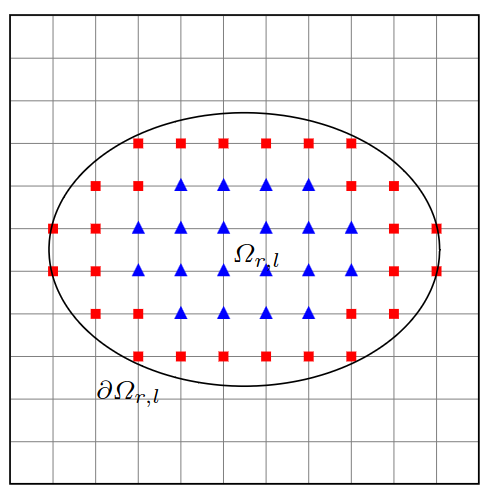}
\caption{A schematic of node identification of the nine-point scheme \eqref{eqn:3d-fully-dis} on a base domain: boundary nodes $\partial\Omega^h_{r,l}$ (marked by red rectangles) and interior nodes $\Omega^h_{r,l}\backslash\partial\Omega^h_{r,l}$(marked by blue triangles).}
\label{fig:grid-nodes}
\end{figure}

\subsection{Boundary condition}
Note that the finite difference schemes \eqref{eqn:2d-fully-dis} and \eqref{eqn:3d-fully-dis} are three-point and nine-point schemes, respectively.
{\color{black}Denote by $\Omega^h_{r,l}$ all grid nodes in $\Omega_{r,l}$.
A grid node is identified as a boundary node if all of its stencil nodes belong to $\Omega_{r,l}^h$.
Otherwise, it is identified as an interior node, see Figure \ref{fig:grid-nodes} for an example.
The set of boundary nodes is denoted as $\partial\Omega^h_{r,l}$.
Further, a control point is identified as a boundary control point if its projection on the base plane belongs to $\partial\Omega_{r,l}^h$.
Denote by $\partial\Gamma_{r,l}^h$ the set of boundary control points.
It is worthwhile to mention that $\partial\Gamma^h_{r,l}$ works as the numerical boundary of $\Gamma_{r,l}$ and is not necessarily a subset of $\partial\Gamma_{r,l}$.}
At boundary nodes, the match condition \eqref{eqn:match-cond} is utilized to pose Dirichlet boundary conditions for \eqref{eqn:2d-pde} and \eqref{eqn:3d-pde}. 
Numerically, the partition of unity $\chi_r$ is taken as
\begin{equation}\label{eqn:pou-num}
    \chi_r(\boldsymbol{x}) = 
    \left \{
    \begin{aligned}
    &1, \quad &\text{if } |\boldsymbol{n}(\boldsymbol{x})\cdot \boldsymbol{e}_r| > |\boldsymbol{n}(\boldsymbol{x})\cdot \boldsymbol{e}_j|, \forall j\neq r, \\
    &0, \quad &\text{otherwise}.
    \end{aligned}
    \right.
\end{equation}
This simple choice of partition of unity is sufficient for accuracy, though it is not a smooth one.
The partition of unity can be understood as evaluating values at boundary nodes on $\Gamma_r$ by interpolation from interior nodes on other subsets $\Gamma_j,j\neq r$. 
Since subsets overlap with each other, interpolation stencils always exist. 

In the following, we introduce the approaches to discretize the matching condition \eqref{eqn:match-cond}.
\subsubsection{Coupled matching condition}
The simplest way to discretize \eqref{eqn:match-cond} is to enforce it at every time level in a discrete sense.
\begin{equation}\label{eqn:cou-match}
    \boldsymbol{x}^{(r), n + 1} = \sum_{j \neq r}\chi_j \boldsymbol{x}^{(j), n+1},\quad\boldsymbol{x}^{(r), n + 1}\in\partial\Gamma_{r}^{h,n+1}.
\end{equation}
This leads to a nonlinear system that couples the solutions on all subsets $\Gamma_r$ in each time step.
Let $\mathbf{u^i}$ and $\mathbf{u^b}$ denote the vectors of solutions at the interior and boundary nodes, respectively. 
The system, which needs to be solved in each time step, is written as
\begin{equation}\label{eqn:sys-1}
    \begin{aligned}
    \mathbf{A}\mathbf{u}^i + \mathbf{Q}\mathbf{u}^b &= \mathbf{f},\\
    \mathbf{u}^b &= \mathbf{\Pi}\mathbf{u}^i,
    \end{aligned}
\end{equation}
where $\mathbf{A}, \mathbf{Q}$  are matrices, $\mathbf{\Pi}$ is the interpolation operator and $\mathbf{f}$ is the vector containing solutions at previous time levels.
Typically, in the system \eqref{eqn:sys-1}, the first equation approximates PDEs, and the second approximates the matching condition.
Here, the operator $\mathbf{\Pi}$ is essentially nonlinear since discretizing \eqref{eqn:cou-match} involves the root-finding of polynomials.
Note that the matrix $\mathbf{A}$ is block-wise diagonal and is invertible.
The nonlinear system \eqref{eqn:sys-1} can be solved, in spirit, with the technique of Schur complement.
One first needs to solve the lower dimensional system
\begin{equation}\label{eqn:schur}
    \mathbf{\Pi}\mathbf{A}^{-1}(\mathbf{f} - \mathbf{Q}\mathbf{u}^b) - \mathbf{u}^b = 0,
\end{equation}
for $\mathbf{u}^b$ and then obtains $\mathbf{u}^i$ by solving
\begin{equation}
    \mathbf{A}\mathbf{u}^i = \mathbf{f} - \mathbf{Q}\mathbf{u}^b.
\end{equation}
The system \eqref{eqn:schur}, which looks like a Schur complement system but is nonlinear, can be solved with the method widely used in domain decomposition methods, the Schwarz alternating method, which is a block-wise Gauss-Seidel type iteration method \cite{mathew2008domain,Saad2003}.
The main idea of the method is to solve problems alternately on each subdomain and to provide boundary conditions for other subdomains.
Generally, the Schwarz alternating method converges geometrically within a few iterations.
The blocks in matrix $\mathbf{A}$ are the approximations of elliptic differential operators, which can also be inverted by an iterative method, such as the successive over-relaxation (SOR) method.
In particular, in two space dimensions, $\mathbf{A}$ is block-wise tri-diagonal, and the Thomas algorithm is applicable.

\subsubsection{ADI method}
Instead of directly enforcing the matching condition \eqref{eqn:match-cond}, we can also follow the idea of the alternating direction implicit (ADI) method, which is used to solve time-dependent PDEs in multiple space dimensions and discretize the matching condition with a time splitting technique.
Note that there is no need to enforce \eqref{eqn:match-cond} accurately since numerical discretization of the PDEs has already introduced numerical errors. 
One only needs to approximate it with an error on the order of $\mathcal{O}(\tau^p)$ where $\tau$ is the time step, and $p$ is the approximation order. 
We evolve $\Gamma_r$ alternately and compute boundary conditions with the newest solutions, which is also an accurate approximation to the matching condition, by
\begin{equation}\label{eqn:adi-match}
    \boldsymbol{x}^{(r), n + 1} = \sum_{j \neq r}\chi_j \boldsymbol{x}^{(j), n^{\star}},\quad \boldsymbol{x}^{(r), n + 1}\in\partial\Gamma_{ r}^{h,n+1}.
\end{equation}
where $\boldsymbol{x}^{(j), n^{\star}}$ is the newest solution on $\Gamma_j$.
For example, boundary conditions for $\Gamma_1$ are interpolated from the newest solution on $\Gamma_2$ and $\Gamma_3$; then, update the solution on $\Gamma_1$ and compute boundary conditions for $\Gamma_2$ using the newest solution on $\Gamma_1$ and $\Gamma_3$, etc.
This approach is non-iterative in the sense that the solutions on subsets $\Gamma_r$ are not coupled, and no Schur complement system needs to be solved in each time step.
This ADI method is also a time-splitting strategy with a formal splitting error on the order of $\mathcal{O}(\tau)$.

\begin{figure}[htbp]
\centering
\subfigure[]{\includegraphics[width=0.3\textwidth]{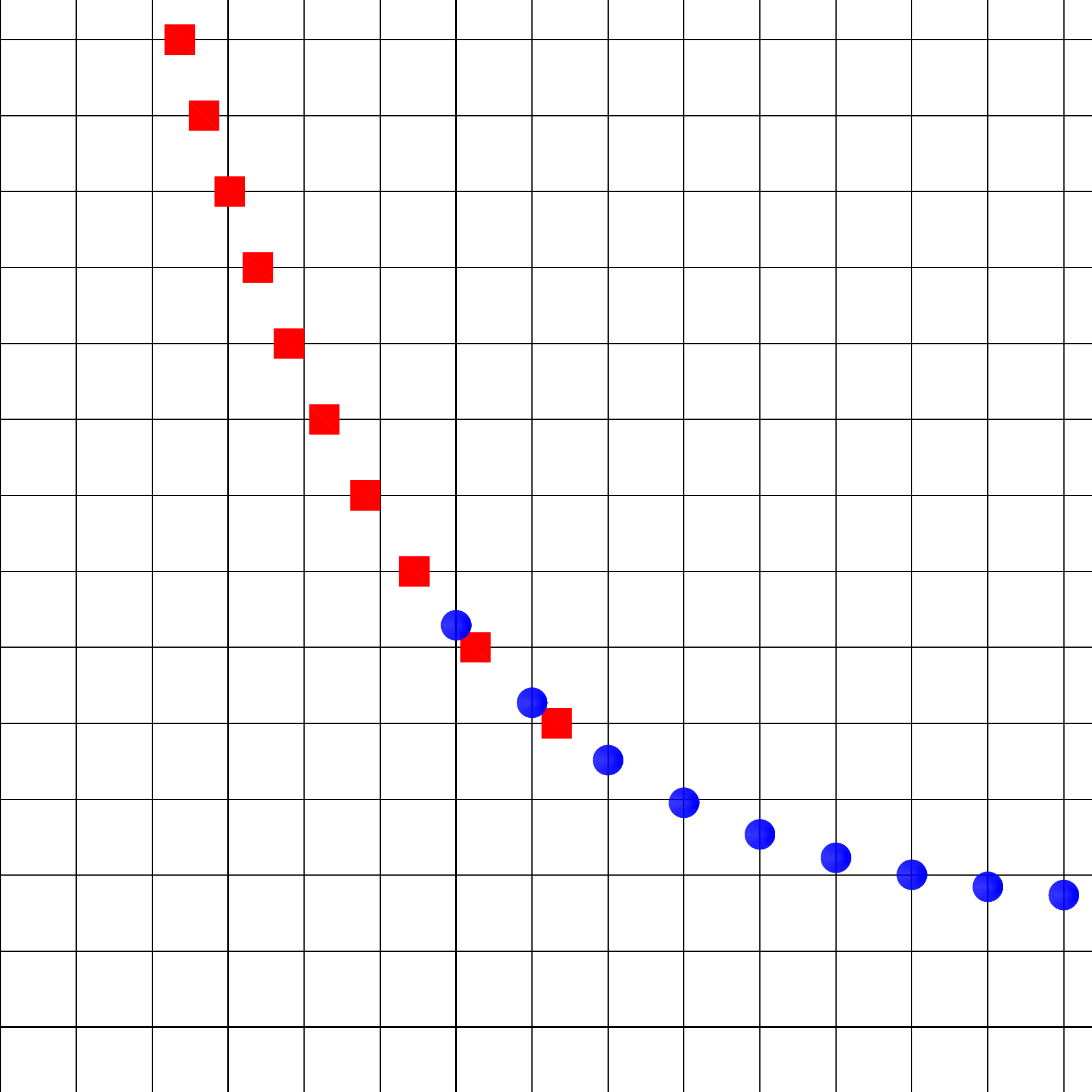}}
\subfigure[]{\includegraphics[width=0.3\textwidth]{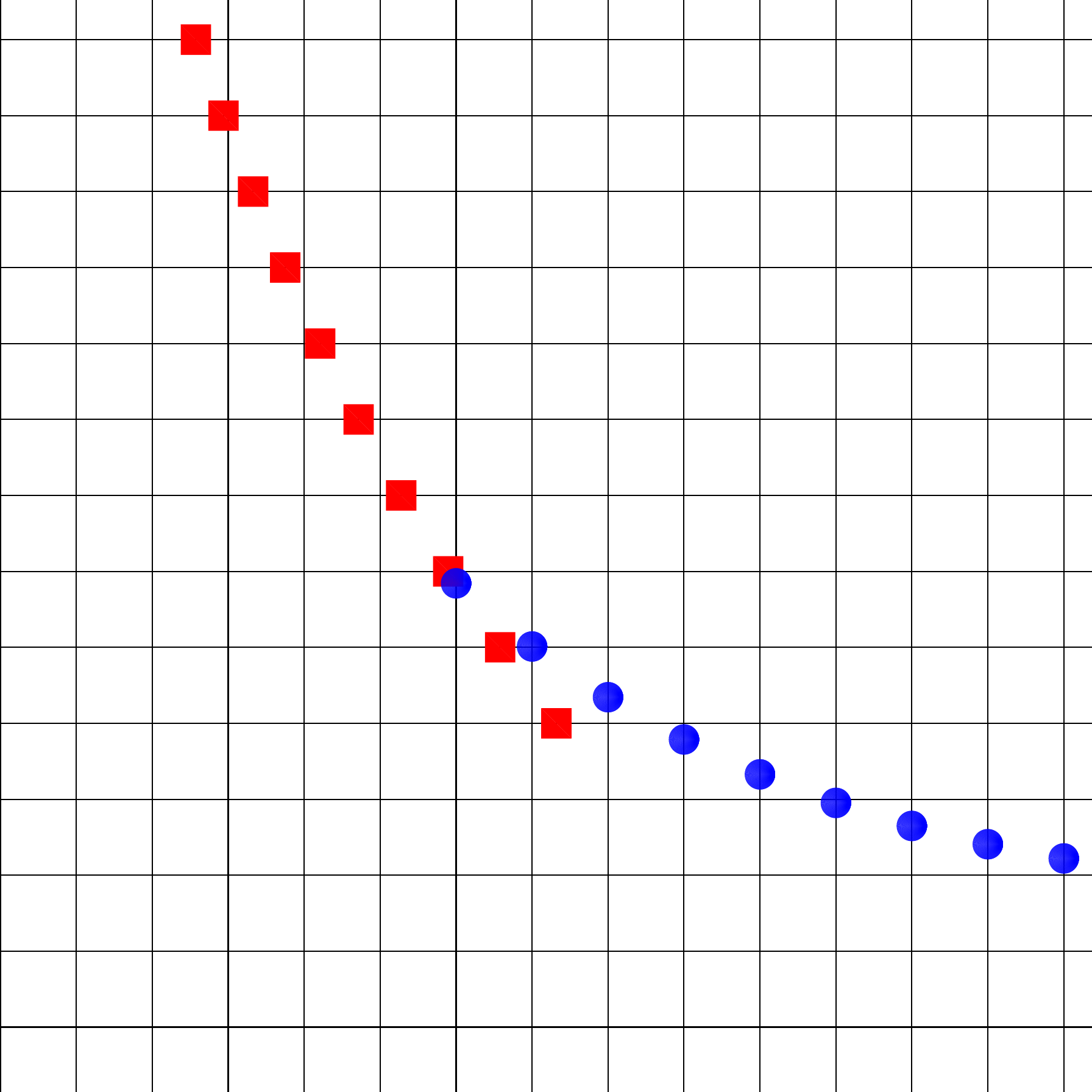}}
\subfigure[]{\includegraphics[width=0.3\textwidth]{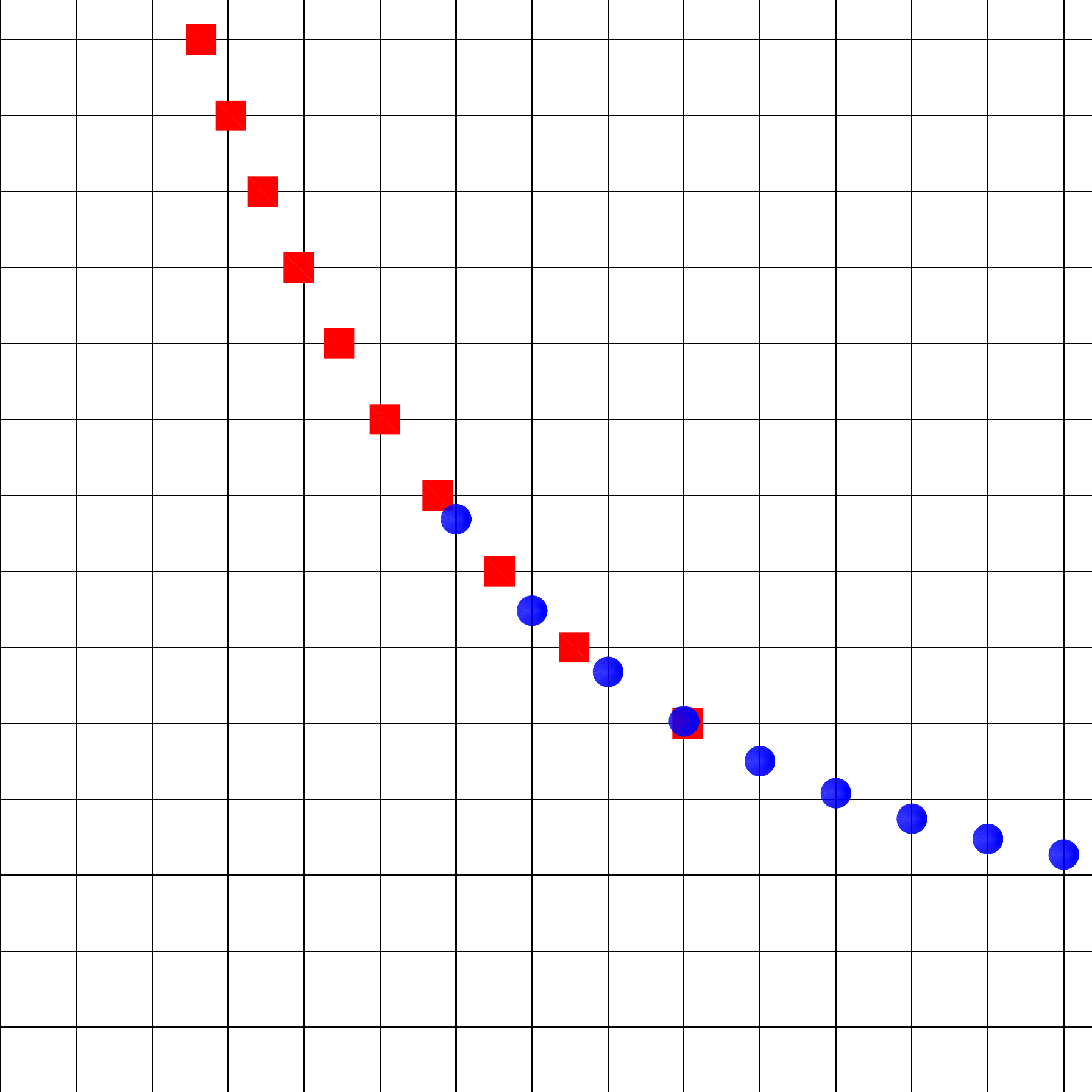}}
\caption{Zoom-in snapshots of numerical solutions. (a) initial solution; (b) solution after a time step by the ADI method; (c) solution after a time step by coupled matching condition. Control points on horizontal lines are marked as red rectangles, and those on vertical lines are marked as blue circles.}
\label{fig:adi-swz}
\end{figure}

The two approaches only differ in the computation of boundary conditions.
Figure \ref{fig:adi-swz} presents the numerical solutions obtained by these two approaches.
One can see that the solutions obtained by these two approaches only have subtle differences at several nodes in the overlapping region. 

In fact, the ADI method and the Schwarz alternating method are closely related to this problem. 
The ADI method is only a strategy to provide Dirichlet boundary conditions for PDEs \eqref{eqn:2d-pde} and \eqref{eqn:3d-pde}.
One can also repeatedly use the ADI method to compute new boundary conditions and update the solution at $t^{n+1}$, which exactly leads to the Schwarz alternating method.
Therefore, the Schwarz alternating method reduces to the ADI method if only one iteration is performed in each time step.

\section{Algorithm summary}\label{sec:algorithm}
In this section, the algorithm for solving mean curvature flows \eqref{eqn:mcf} with the proposed ADI method is summarized as follows:
\begin{description}
\item[Algorithm.] The ADI method for mean curvature flows:
\item[Step 1.] Given the initial hypersurface by its parametric form or a level set function, embed it into a bounding box that is uniformly partitioned into a Cartesian grid and find the control points with overlapping decomposition strategy described in Subsection \ref{sec:osd-dis}. 
\item[Step 2.] In each time step, evolve the overlapping subsets alternately. For each subset, do the following three procedures:
\begin{enumerate}
    \item[(a)] identify all the nodes in each isolated component of the subset with the breadth-first search method.
    \item[(b)] compute the Dirichlet-type boundary condition for boundary nodes with the discrete matching condition \eqref{eqn:adi-match};
    \item[(c)] evolve the subset to the next time level by solving \eqref{eqn:2d-fully-dis} or \eqref{eqn:3d-fully-dis};
\end{enumerate}
\item[Step 4.] Update control points such that they satisfy the decomposition strategy described in Subsection \ref{sec:osd-dis}.
\item[Step 5.] Repeat steps 2-4 until the final computational time.
\end{description}
\begin{remark}
Procedure (a) in Step 2 is only for matrix assembly such that direct methods, such as the Thomas algorithm, are applicable.
Suppose the finite difference equations \eqref{eqn:2d-fully-dis} and \eqref{eqn:3d-fully-dis} are solved with iterative methods which only require the matrix-vector product. In that case, one can find the local stencil points on-the-fly instead of finding all the nodes in the components in advance.
\end{remark}

\section{Numerical results}\label{sec:result}
This section presents numerical examples in two and three space dimensions to validate the proposed method. 
Initial hypersurfaces are given in parametric forms or the zero level set of level set functions, which will be prescribed in each example. 
In all the numerical examples, the bounding box $\mathcal{B}$ is uniformly partitioned into a Cartesian grid $\mathcal{G}$ with $N$ intervals in each direction.
For problems with exact solutions, we estimate the numerical error at a surface point by finding its projection on the hypersurface of the exact solution and computing the distance.
We take the solution on a fine grid for problems without exact solutions as a reference ``exact'' solution.
Then, we estimate the numerical error at a surface point by finding the closest surface point on the reference solution and computing the distance.
The numerical errors in the maximum norm and $l_2$ norm are computed by
\begin{equation}
    \Vert \mathbf{e}_h\Vert_{\infty} = \max_{\boldsymbol{x}_i\in\Gamma^h}\left\{\Vert \boldsymbol{x}_i - \boldsymbol{x}^{ref}_i\Vert\right\}, \quad \Vert \mathbf{e}_h\Vert_{2} = \sqrt{\dfrac{1}{N_{\Gamma}} \sum_{\boldsymbol{x}_i\in\Gamma^h}\Vert \boldsymbol{x}_i - \boldsymbol{x}^{ref}_i\Vert^2 },
\end{equation}
where $\boldsymbol{x}^{ref}_i$ is the exact solution or reference solution on a fine grid associated with $\boldsymbol{x}_i$ and $N_{\Gamma}$ is the total point number.

The following numerical experiments are performed on a personal computer with a 3.80 GHz Intel Core i7 processor. 
The codes for conducting the numerical experiments are written in C++ computer language.

\subsection{Two space dimensional examples}\label{eg:2d-egs}
First, we solve the mean curvature flow for a simple case and compare the numerical solution with the exact solution to verify the convergence of the proposed method.
The initial shape is chosen such that the curve is a circle whose radius $r(t)$ satisfies
\begin{equation}\label{eqn:init-cir}
r(t) = \sqrt{1 - 2t}.
\end{equation}
The bounding box is taken as $[-1.2, 1.2]^2$. 
Note that the curve will eventually shrink to a point at $T_{end} = 0.5$.
We chose to estimate the numerical errors at $T = 0.2$ to ensure that the coarsest grid $N = 64$ can fully resolve the curve during the computation.
On finer grids, the computation can last longer than $T$.
Time step size is chosen to be $\Delta t = 0.1\Delta x$ where $\Delta x$ is the spatial grid size. 
Numerical results are summarized in Table \ref{tab:cir-err}.

\begin{table}[htbp]
\caption{Numerical error and convergence order of the 2D MCF for a circle-shaped initial curve.}
\label{tab:cir-err}
\centering
\begin{tabular}{|c|c|c|c|c|}
\hline
N   & $\Vert \boldsymbol{e}_h\Vert_{\infty}$ & order  &  $\Vert \boldsymbol{e}_h\ \Vert_{2}$ &  order \\ \hline
64  & 7.99e-03 & - & 3.99e-03 & - \\ \hline
128  & 3.80e-03 & 1.07 & 1.79e-04 & 1.16 \\ \hline
256  & 1.88e-03 & 1.02 & 8.42e-04 &  1.09 \\ \hline
512  & 9.45e-04 &  0.99 & 4.12e-05 & 1.03 \\ \hline
1024 & 4.72e-04 &  1.00 & 2.05e-05 &  1.01 \\ \hline
\end{tabular}
\end{table}

\begin{figure}[htbp]
\centering
\subfigure[t = 0 ]{\includegraphics[scale=0.1]{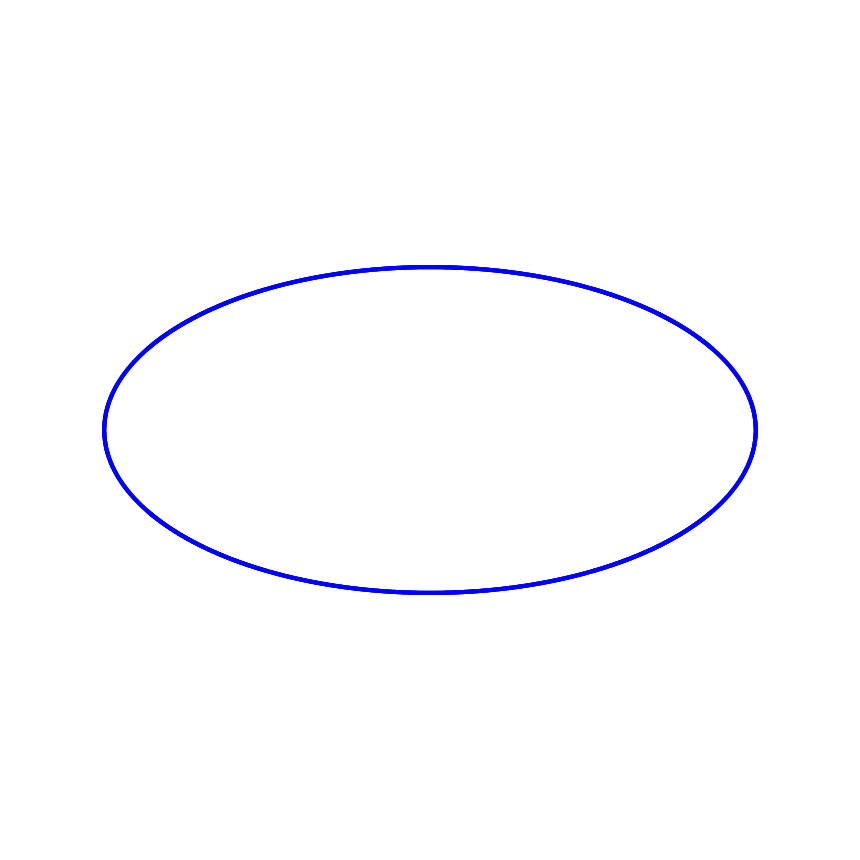}}
\subfigure[t = 0.04 ]{\includegraphics[scale=0.1]{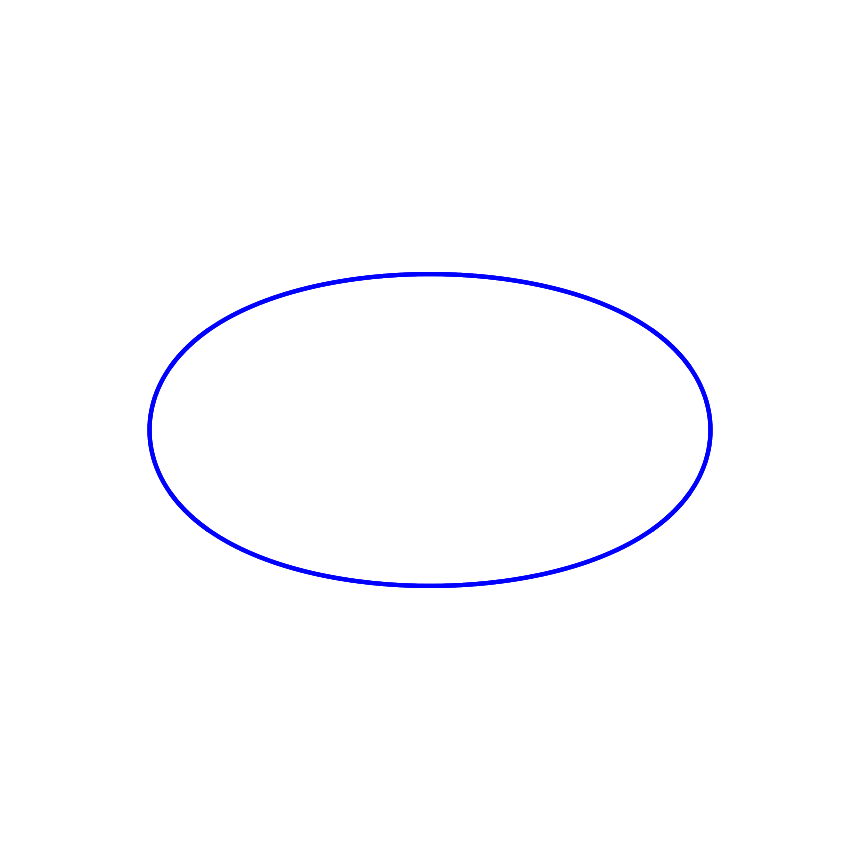}}
\subfigure[t = 0.08 ]{\includegraphics[scale=0.1]{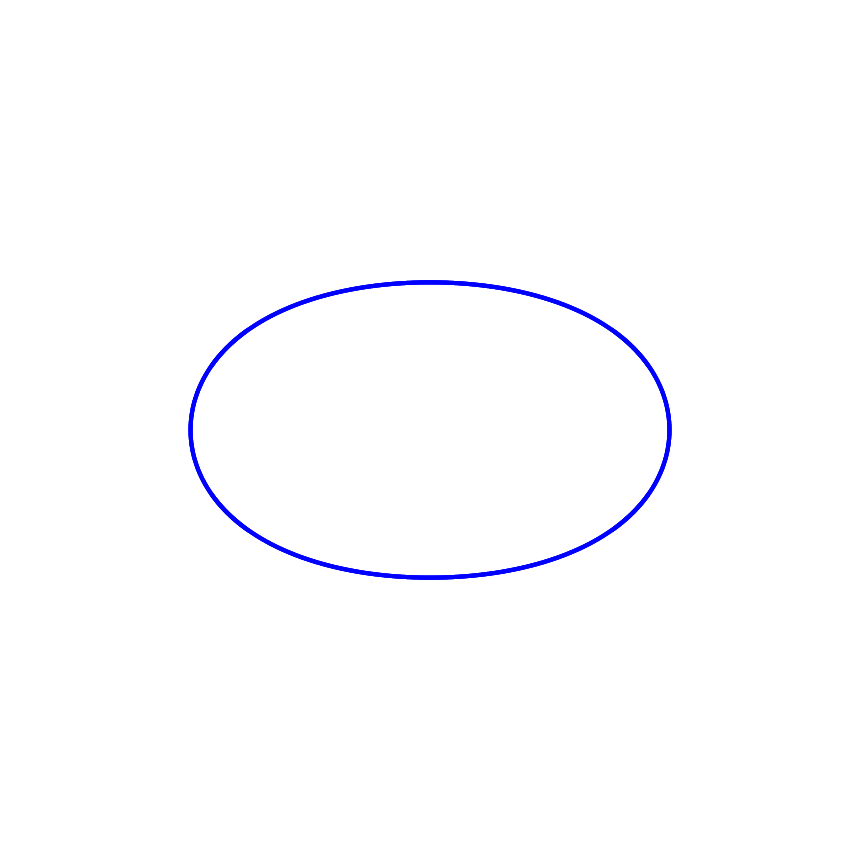}}\\
\subfigure[t = 0.12 ]{\includegraphics[scale=0.1]{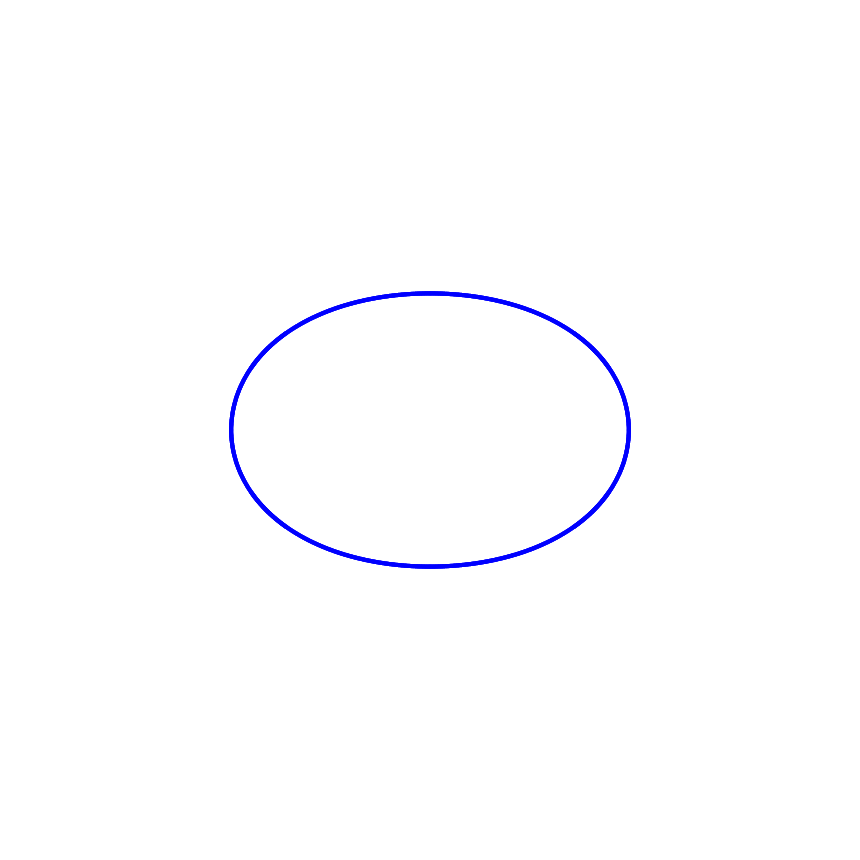}}
\subfigure[t = 0.16 ]{\includegraphics[scale=0.1]{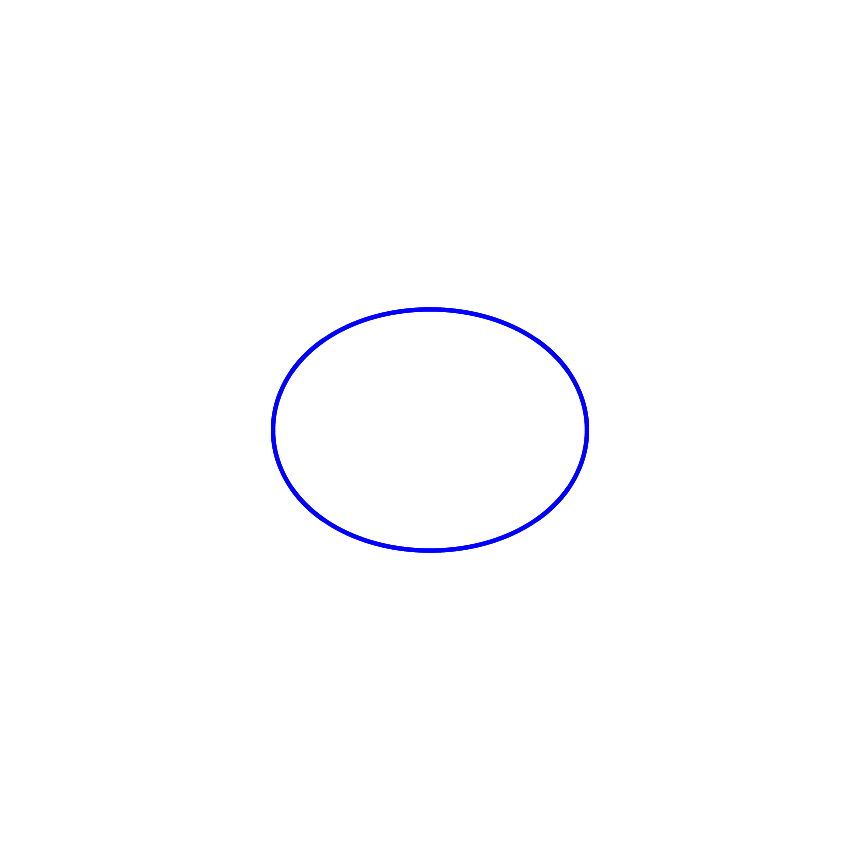}}
\subfigure[t = 0.20 ]{\includegraphics[scale=0.1]{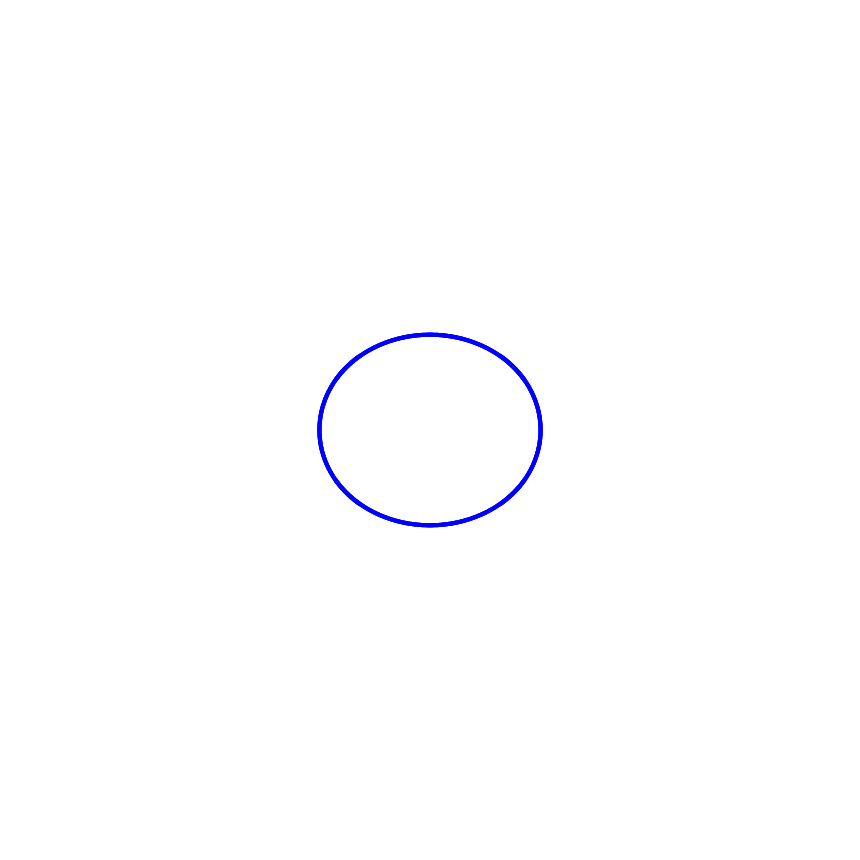}}
\caption{Time evolution of the MCF for an ellipse-shaped curve.}
\label{fig:ellipse}
\end{figure}

\begin{figure}[htbp]
    \centering
    \includegraphics[width=0.8\textwidth]{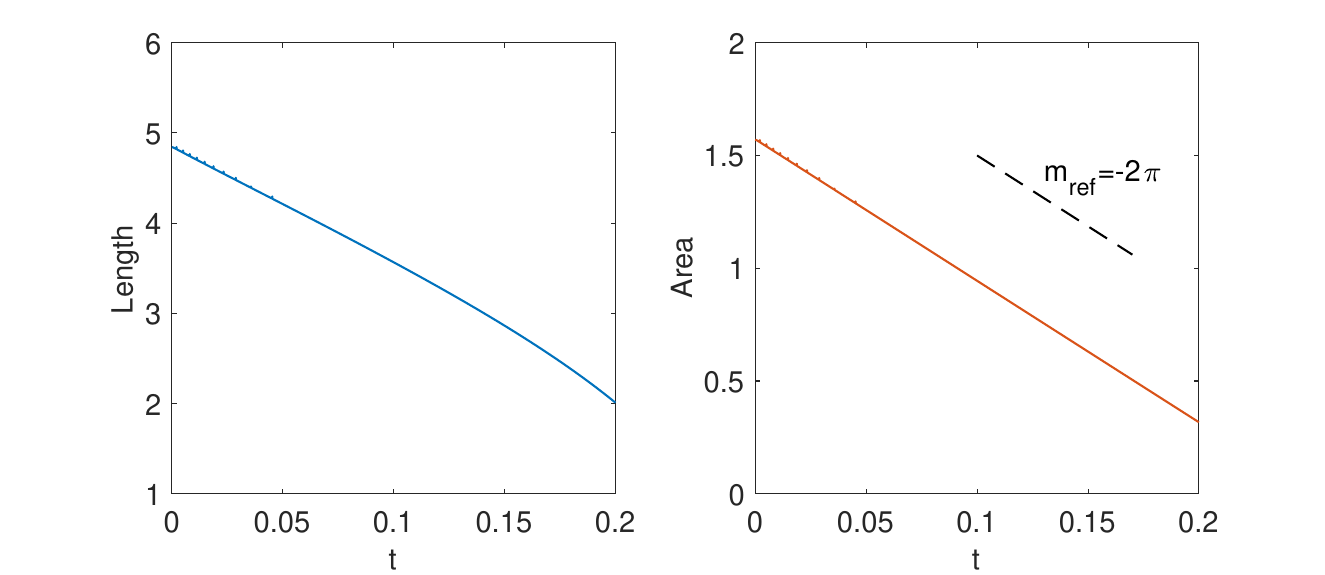}
    \caption{Time evolution of curve length and enclosed area of the ellipse-shaped curve.}
    \label{fig:ell-L-A}
\end{figure}

Next, we change the initial shape to an ellipse which is given by
\begin{equation}\label{eqn:init-ellipse}
\left \{
\begin{aligned}
x &= a \cos(\theta),\\
y &= b \sin(\theta),
\end{aligned}
\right.
\quad \theta \in[0, 2\pi),
\end{equation}
with $a = 1.0, b = 0.5$. 
The problem is solved in the bounding box $\mathcal{B}=[-1.2,1.2]^2$.
Time step size is taken as $\Delta t = 0.1\Delta x$.
Since there is no exact solution for this configuration, the solution on a fine grid with $N = 2048$ is chosen as a reference solution to estimate numerical errors. 
The estimated error and the convergence order are summarized in Table \ref{tab:ellipse}.
It can be observed that the convergence order is a bit larger than $1$.
This may be due to the inaccurate estimation of numerical error based on the distance between the surface point to its closest point in the reference solution.
The evolution history of the curve is presented in Figure \ref{fig:ellipse}.
The changes in curve length and enclosed area are computed on the grid with $N=1024$ and shown in Figure \ref{fig:ell-L-A}.
It can be observed that the enclosed area loss rate compares favorably with the theoretic result, whose slope is $m_{ref} = -2\pi$.

\begin{table}[htbp]
\caption{Numerical error and convergence order of the 2D MCF for an ellipse-shaped initial curve.}
\label{tab:ellipse}
\centering
\begin{tabular}{|c|c|c|c|c|c|}
\hline
$N$     & $\Vert \boldsymbol{e}_h\Vert_{\infty}$ & order  &  $\Vert \boldsymbol{e}_h \Vert_{2}$ &  order \\ \hline
128  & 2.45e-02 & - & 1.84e-02 & - \\ \hline
256  & 1.04e-02 & 1.24 & 7.95e-03 &  1.21\\ \hline
512  & 4.31e-03 &  1.27 & 3.31e-03 & 1.26\\ \hline
1024 & 1.55e-03 &  1.48 & 1.15e-03 &  1.52\\ \hline
\end{tabular}
\end{table}


\begin{figure}[htbp]
\centering
\subfigure[t = 0 ]{\includegraphics[scale=0.1]{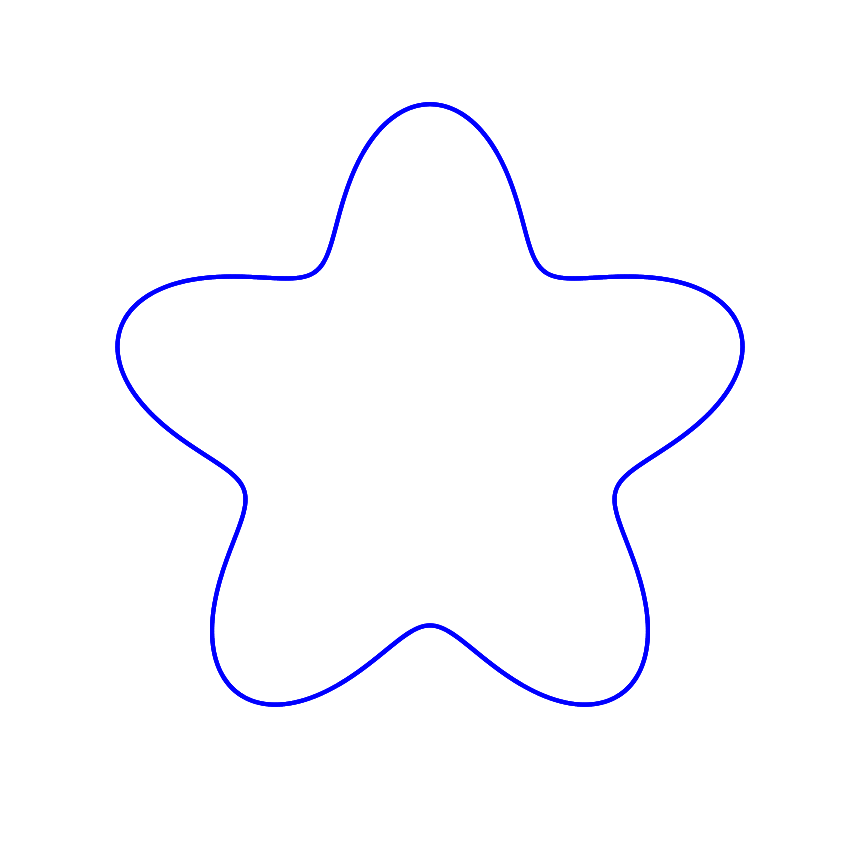}}
\subfigure[t = 0.02 ]{\includegraphics[scale=0.1]{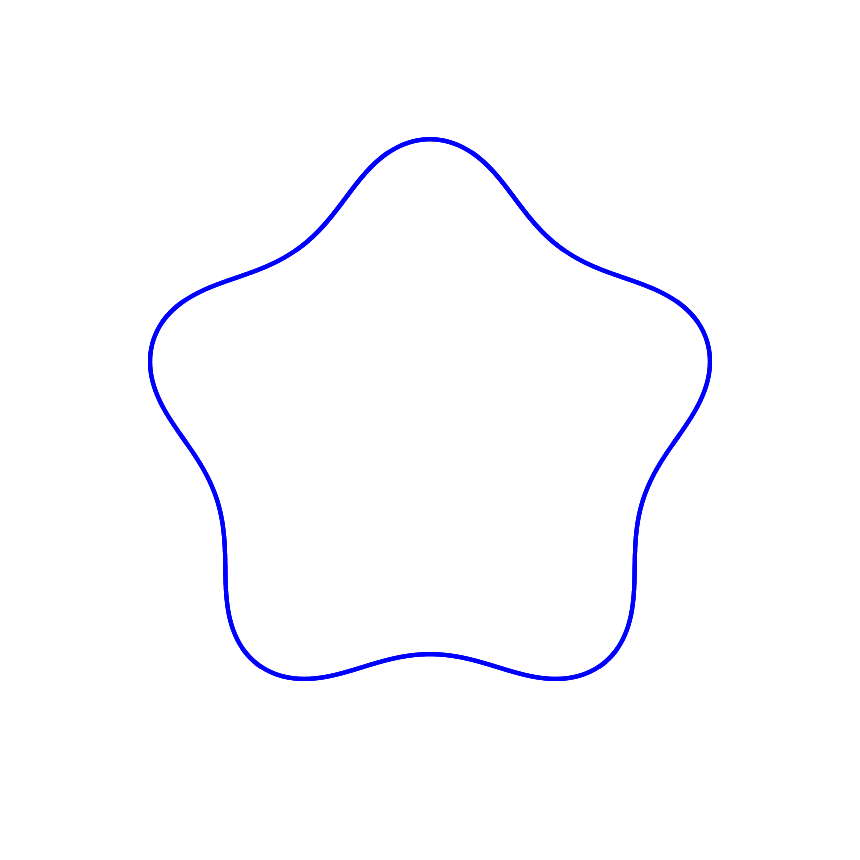}}
\subfigure[t = 0.04 ]{\includegraphics[scale=0.1]{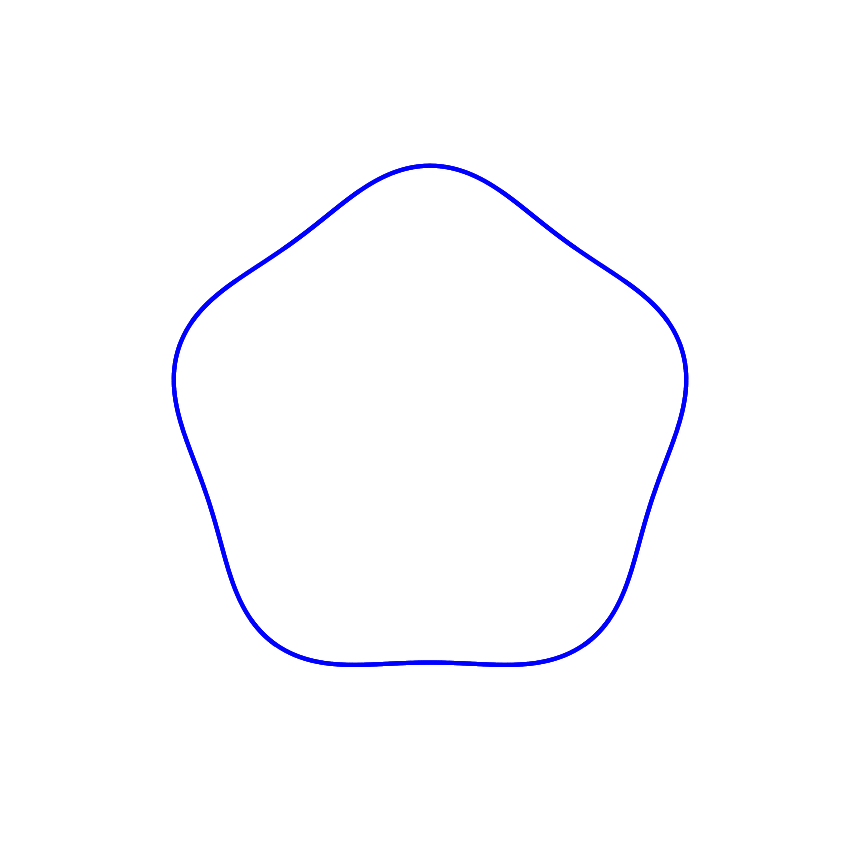}}\\
\subfigure[t = 0.06 ]{\includegraphics[scale=0.1]{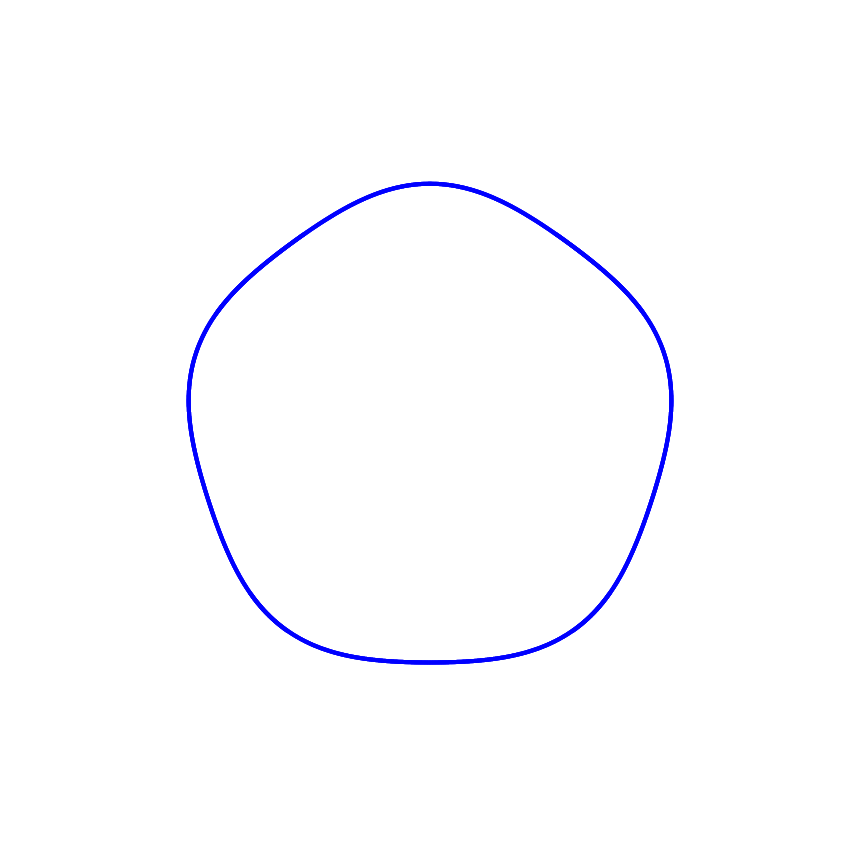}}
\subfigure[t = 0.08 ]{\includegraphics[scale=0.1]{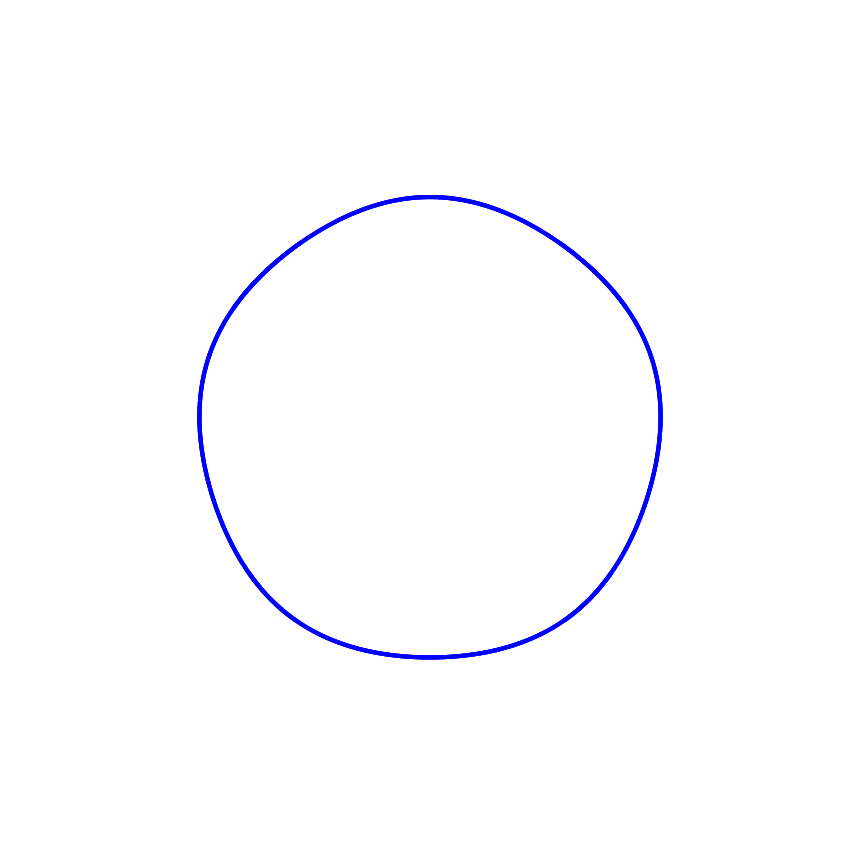}}
\subfigure[t = 0.10 ]{\includegraphics[scale=0.1]{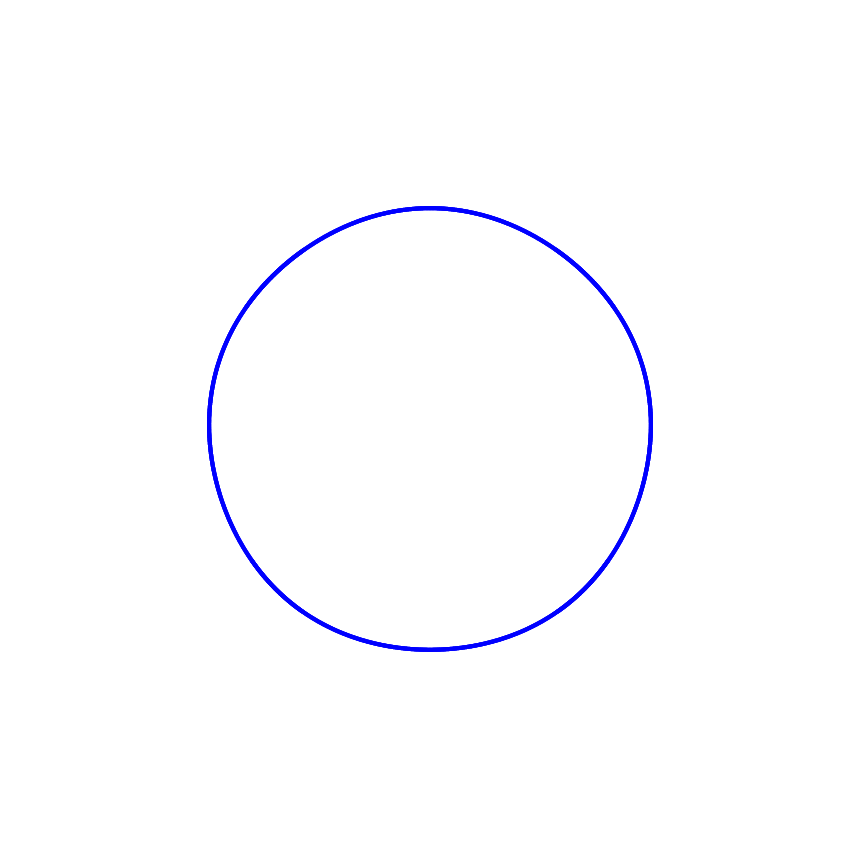}}
\caption{Time evolution of the MCF for a five-fold star-shaped curve.}
\label{fig:star}
\end{figure}

We also chose a five-fold star-shaped initial curve, which is given by
\begin{equation}\label{eqn:init-star}
\left \{
\begin{aligned}
x &= a(\kappa + \eta \sin(m\theta))\cos(\theta),\\
y &= b(\kappa + \eta \sin(m\theta))\sin(\theta),
\end{aligned}
\right.
\quad \theta \in[0, 2\pi),
\end{equation}
with $a = 1.0, b = 1.0, \kappa = 0.8, \eta = 0.2$ and $m = 5$. 
The bounding box and time step size are the same as those in the last case.
Estimated numerical errors and convergence orders are summarized in Table \ref{tab:star}.
The evolution history of the curve and changes in curve length and enclosed area are presented in Figure \ref{fig:star} and \ref{fig:ell-L-A}, respectively.
The numerical results are also consistent with theoretical results.

\begin{table}[htbp]
\caption{Numerical error and convergence order of the 2D MCF for a five-fold star-shaped initial curve.}
\label{tab:star}
\centering
\begin{tabular}{|c|c|c|c|c|c|}
\hline
$N$    & $\Vert \boldsymbol{e}_h\Vert_{\infty}$ & order  &  $\Vert \boldsymbol{e}_h \Vert_{2}$ &  order \\ \hline
128  & 5.45e-03 & - & 3.21e-03 & -  \\ \hline
256  & 2.60e-03 & 1.07 & 1.46e-03 &  1.14 \\ \hline
512  & 1.21e-03 &  1.10 & 6.60e-04 & 1.15 \\ \hline
1024 & 5.00e-04 &  1.28 & 2.26e-05 &  1.55 \\ \hline
\end{tabular}
\end{table}

\begin{figure}[htbp]
    \centering
    \includegraphics[width=0.8\textwidth]{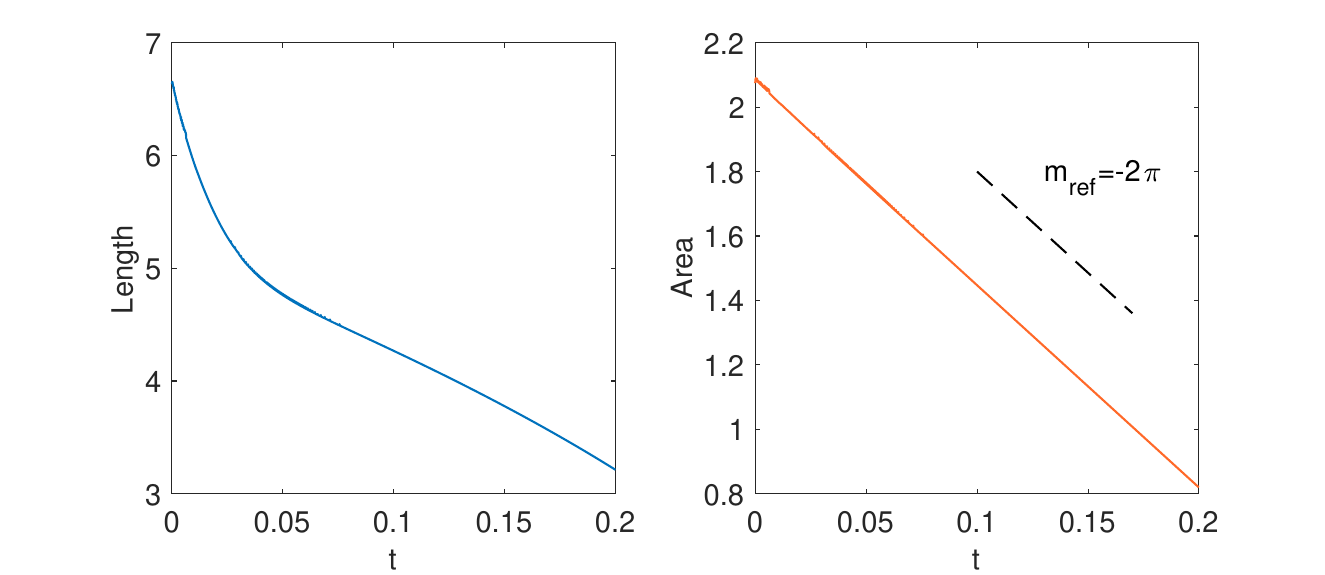}
    \caption{Time evolution of curve length and enclosed area of the five-fold star-shaped curve.}
    \label{fig:star-L-A}
\end{figure}
For this case, we compare the present method with an explicit time-advancing scheme for the mean curvature flow, which discretizes the equation
\begin{equation}
\frac{d}{dt}
\begin{pmatrix}
x\\
y
\end{pmatrix}
= - \frac{x_{\theta}y_{\theta\theta} - x_{\theta\theta}y_{\theta}}{(x_{\theta}^{2} + y_{\theta}^{2})^{2}}
\begin{pmatrix}
y_{\theta} \\
- x _{\theta}
\end{pmatrix},\quad \theta\in[0, 2\pi),
\end{equation}
with forward Euler scheme and central differences for temporal and spatial derivatives, respectively.
The time step size for the explicit method is chosen to ensure numerical stability, using the adaptive time step 
\begin{equation}
\Delta t = 0.8(\min_{\theta}\Delta s)^2,
\end{equation}
where $\Delta s$ denotes the Euclidean distance between two adjacent grid nodes.
We properly optimize the codes for both methods and collect the required CPU times for solving the mean curvature flow to the final time $T = 0.1$.
The control point number at time $t$ is denoted by $M_t$.
Numerical results are summarized in Table \ref{tab:cpu-time}.
It can be observed that while the forward Euler method is faster on coarse grids, the present method becomes more efficient than the forward Euler method as the point number increases.
It can be explained by the complexities of the two methods.
To compute the solution to a fixed final time, since the time step size can be chosen as linearly proportional to the spatial grid size for the present method, the computational complexity is $\mathcal{O}(N^2)$. 
However, the computational complexity of the forward Euler scheme is $\mathcal{O}(N^3)$ due to high-order constraints on time step size.
In fact, the forward Euler method fails for long-time computation since the required time step size quickly decreases to $10^{-6}$ due to the curve shortening phenomenon and poor mesh quality.

\begin{table}[htbp]
\centering
\caption{CPU time comparison between the present method and an explicit time advancing scheme.}
\label{tab:cpu-time}
\begin{tabular}{|c|c|c|c|c|c|}
\hline
\multicolumn{4}{|c|}{Present method}                             & \multicolumn{2}{c|}{Forward Euler method} \\ \hline
$N$ & $M_0$ & $M_T$ & CPU  times(secs) & $M_0=M_T$ & CPU times(secs)       \\ \hline
128  & 392  & 224  & 9.66e-03 & 224  & 3.73e-03 \\ \hline
256  & 786  & 444  & 3.54e-02 & 448  & 1.73e-02  \\ \hline
512  & 1574 & 890  & 1.34e-01 & 896  & 1.34e-01 \\ \hline
1024 & 3150 & 1776 & 5.15e-01 & 1792 & 1.02e+00 \\ \hline
2048 & 6300 & 3548 & 2.08e+00 & 3584 & 8.30e+00 \\ \hline
\end{tabular}
\end{table}

\subsection{Three space dimensional examples}\label{eg:3d-egs}
For three space dimensional mean curvature flows, we test the convergence rate of the method by considering a simple case, a sphere-shaped initial surface.
Similar to two space dimensional case, this configuration has an exact solution: the surface maintains a sphere, and the radius $r(t)$ satisfies
\begin{equation}
r(t) = \sqrt{1 - 2t}.
\end{equation}
Numerical errors are estimated at $T=0.2$ and summarized in Table \ref{tab:sph-err}.
\begin{table}[htbp]
\caption{Numerical error and convergence order of the 3D MCF for a sphere-shaped initial surface.}
\label{tab:sph-err}
\centering
\begin{tabular}{|c|c|c|c|c|}
\hline
N   & $\Vert \boldsymbol{e}_h\Vert_{\infty}$ & order  &  $\Vert \boldsymbol{e}_h \Vert_{2}$ &  order \\ \hline
64  & 4.72e-03 & - & 1.42e-03 & - \\ \hline
128  & 3.46e-03 & 0.45 & 7.46e-04 & 0.93 \\ \hline
256  & 1.45e-03 & 1.25 & 3.63e-04 &  1.04 \\ \hline
512  & 9.18e-04 & 0.66  & 1.82e-04 &  1.00 \\ \hline
1024 & 3.47e-04 &  1.40 & 8.81e-05 &  1.05 \\ \hline
\end{tabular}
\end{table}

We also solve the mean curvature flow in three space dimensions for more examples.
In the following numerical examples, the bounding boxes partitioned into a Cartesian grid are all chosen as $\mathcal{B} = [-1.2, 1.2]^3$.
The time step is chosen as $\Delta t = 0.05\Delta x$.

\begin{figure}[htbp]
\centering
\subfigure[t = 0 ]{\includegraphics[width=0.3\textwidth]{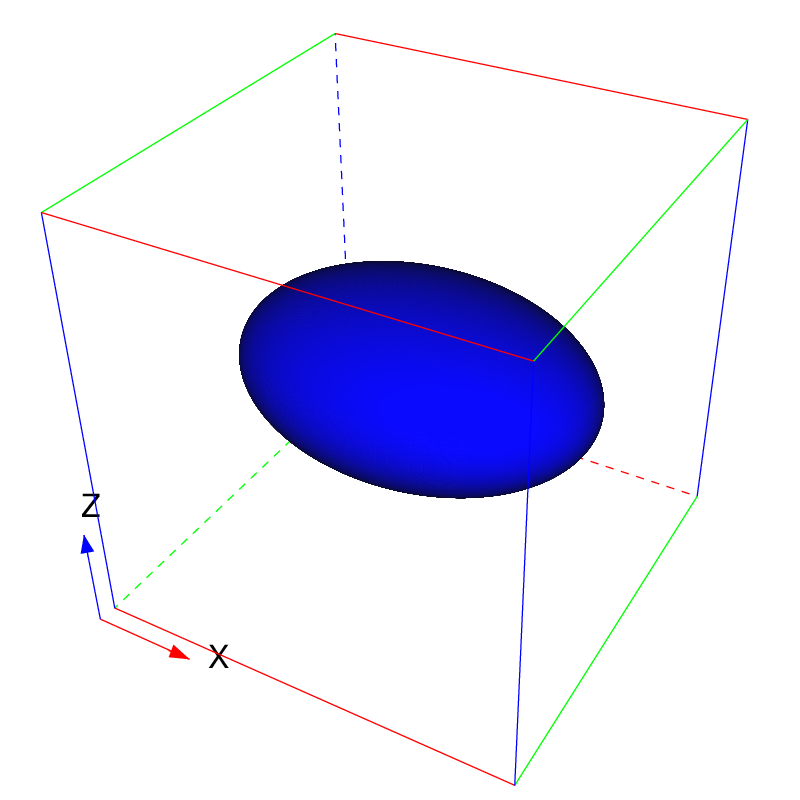}}
\subfigure[t = 0.04 ]{\includegraphics[width=0.3\textwidth]{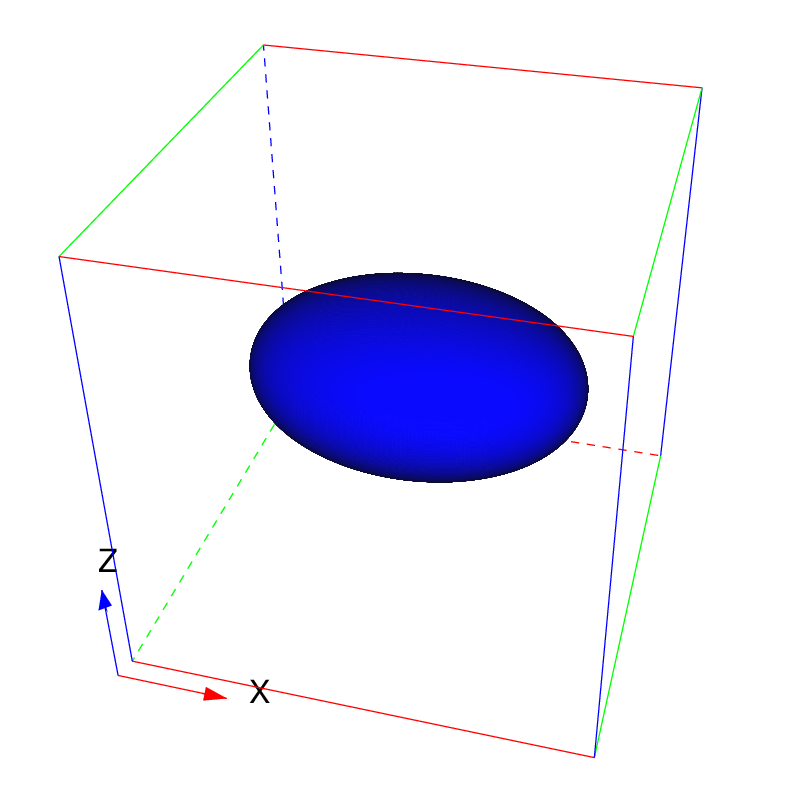}}
\subfigure[t = 0.08 ]{\includegraphics[width=0.3\textwidth]{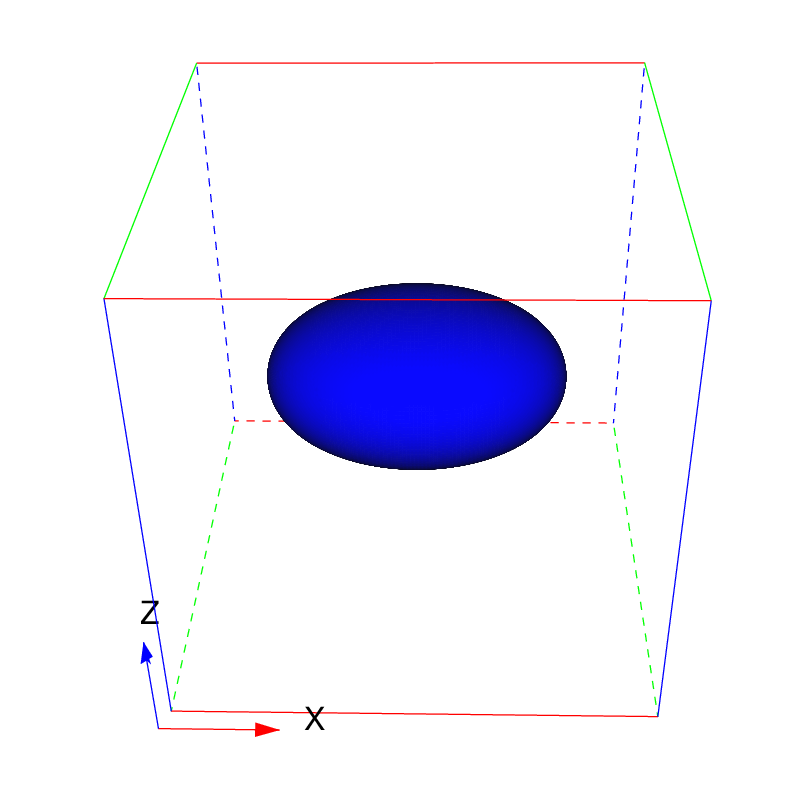}}\\
\subfigure[t = 0.12 ]{\includegraphics[width=0.3\textwidth]{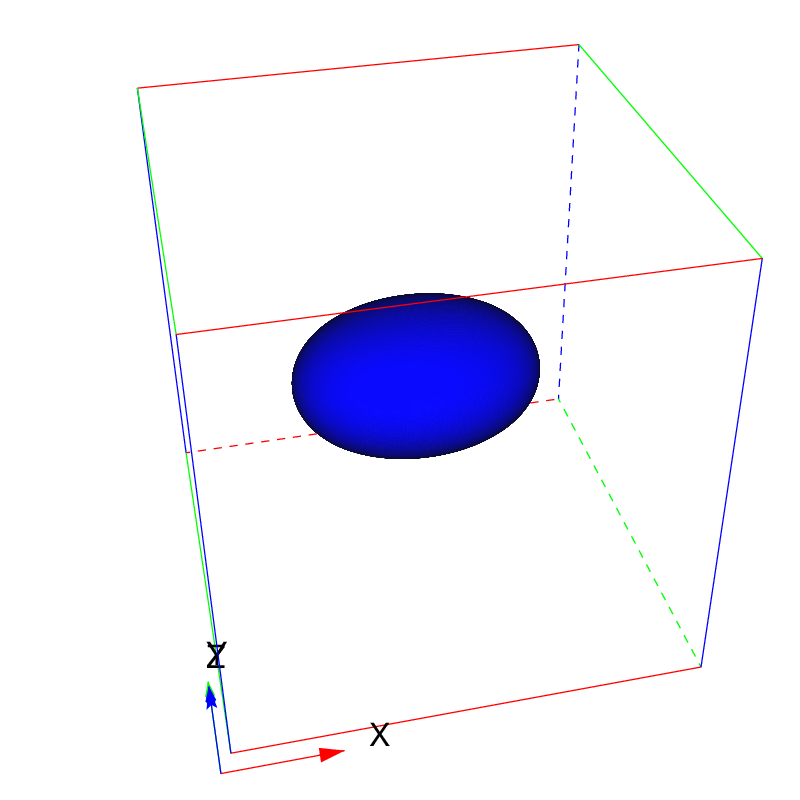}}
\subfigure[t = 0.16 ]{\includegraphics[width=0.3\textwidth]{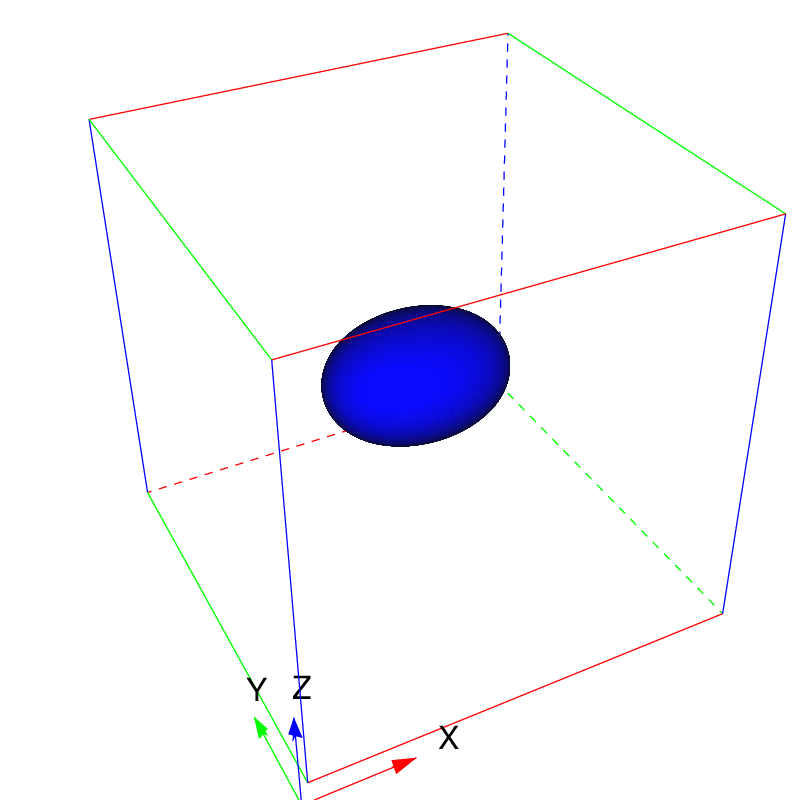}}
\subfigure[t = 0.20]{\includegraphics[width=0.3\textwidth]{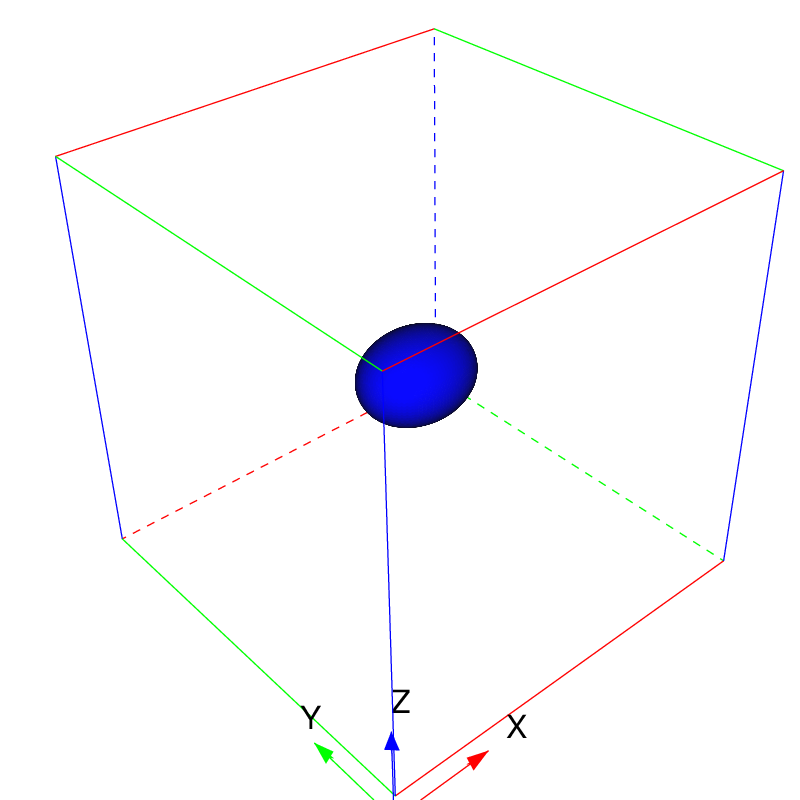}}
\caption{Time evolution of the MCF for an ellipsoid-shaped surface.}
\label{fig:ellipsoid}
\end{figure}

\begin{figure}[htbp]
    \centering
    \includegraphics[width=0.8\textwidth]{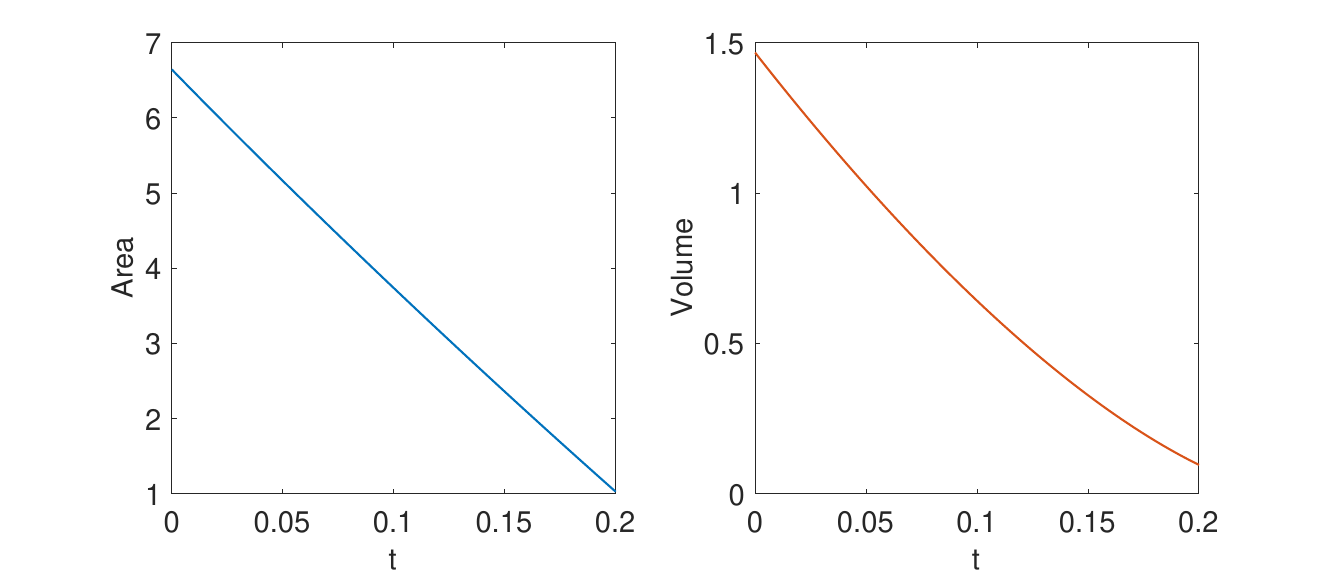}
    \caption{Time evolution of the surface area and enclosed volume of the ellipsoid-shaped surface.}
    \label{fig:ell-A-V}
\end{figure}

In the first case, we set the initial shape as an ellipsoid which is given by
\begin{equation}
\Gamma = \left\{(x,y,z)\Big |\frac{x^2}{a^2} + \frac{y^2}{b^2} + \frac{z^2}{c^2} - 1 = 0\right\},
\end{equation}
with $a = 1.0, b = 0.7, c = 0.5$.
Numerical error and convergence order estimated at $t = 0.2$ are summarized in Table \ref{tab:ellipsoid}.
The time evolution of the surface and its area and enclosed volume are presented in Figure  \ref{fig:ellipsoid} and \ref{fig:ell-A-V}, respectively.
One can observe that the major axis of the ellipsoid decreases faster compared with the other two axes and the ellipsoid becomes very close to a sphere.
Surface area decreases with time, which is consistent with the theoretical result. 

\begin{table}[htbp]
\caption{Numerical error and convergence order of the 3D MCF for an ellipsoid-shaped initial surface.}
\label{tab:ellipsoid}
\centering
\begin{tabular}{|c|c|c|c|c|}
\hline
$N$     & $\Vert \boldsymbol{e}_h\Vert_{\infty}$ & order  &  $\Vert \boldsymbol{e}_h \Vert_{2}$ &  order \\ \hline
128  & 1.58e-03 & - & 3.78e-04 & - \\ \hline
256  & 1.32e-03 & 0.26 & 1.74e-04 &  1.12 \\ \hline
512  & 3.70e-04 & 1.83  & 6.36e-05 &  1.45 \\ \hline
\end{tabular}
\end{table}

\begin{figure}[htbp]
\centering
\subfigure[t = 0 ]{\includegraphics[width=0.3\textwidth]{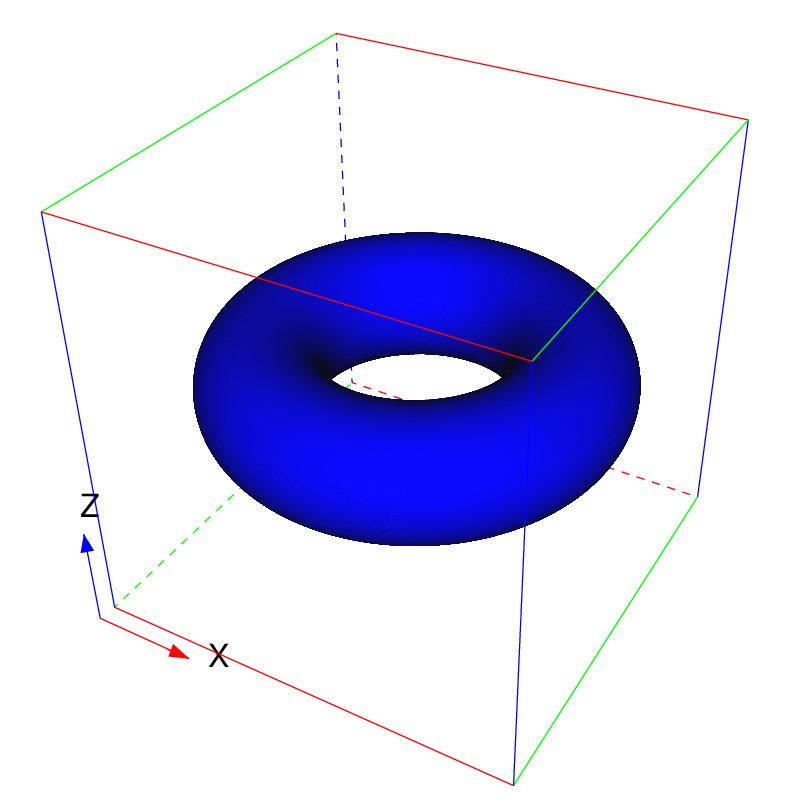}}
\subfigure[t = 0.02 ]{\includegraphics[width=0.3\textwidth]{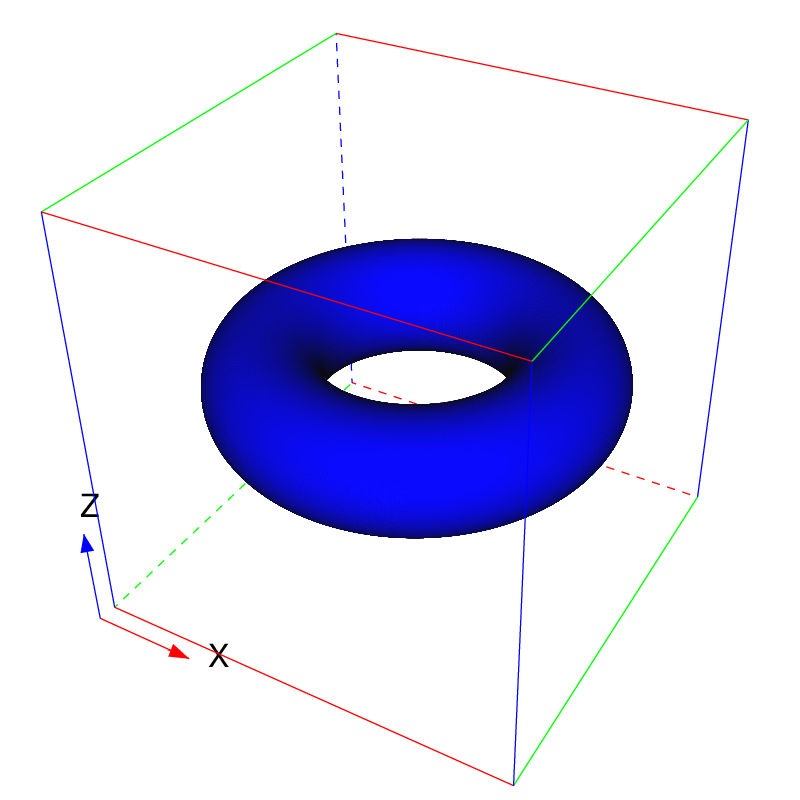}}
\subfigure[t = 0.04 ]{\includegraphics[width=0.3\textwidth]{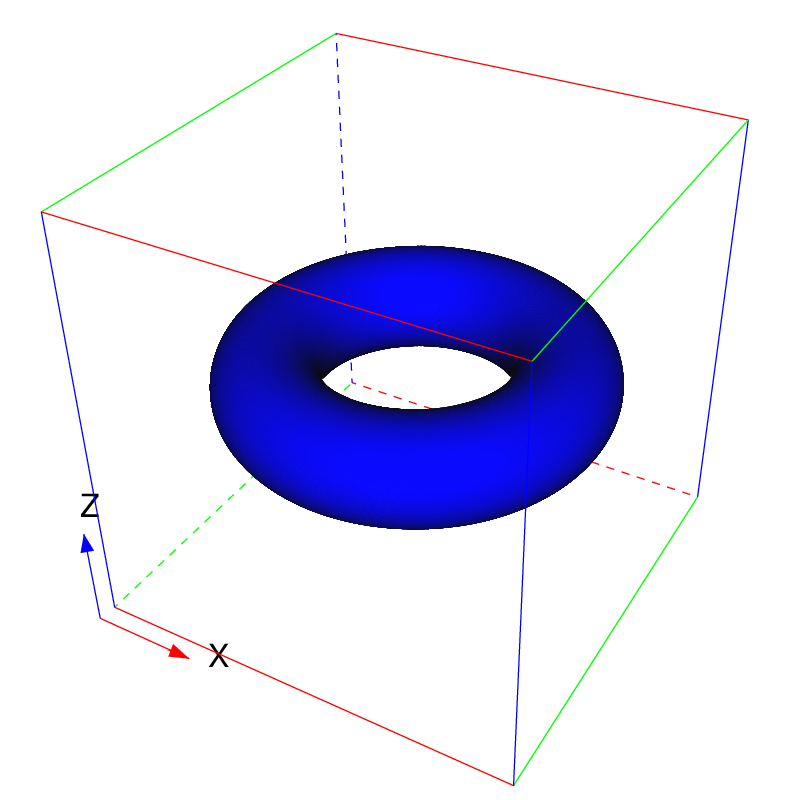}}
\subfigure[t = 0.06 ]{\includegraphics[width=0.3\textwidth]{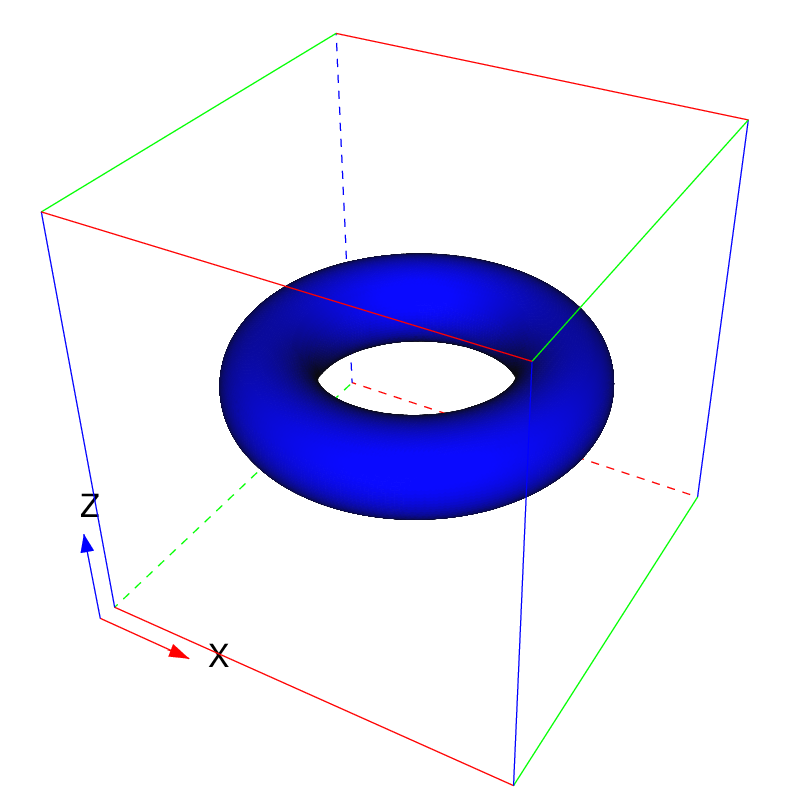}}
\subfigure[t = 0.08 ]{\includegraphics[width=0.3\textwidth]{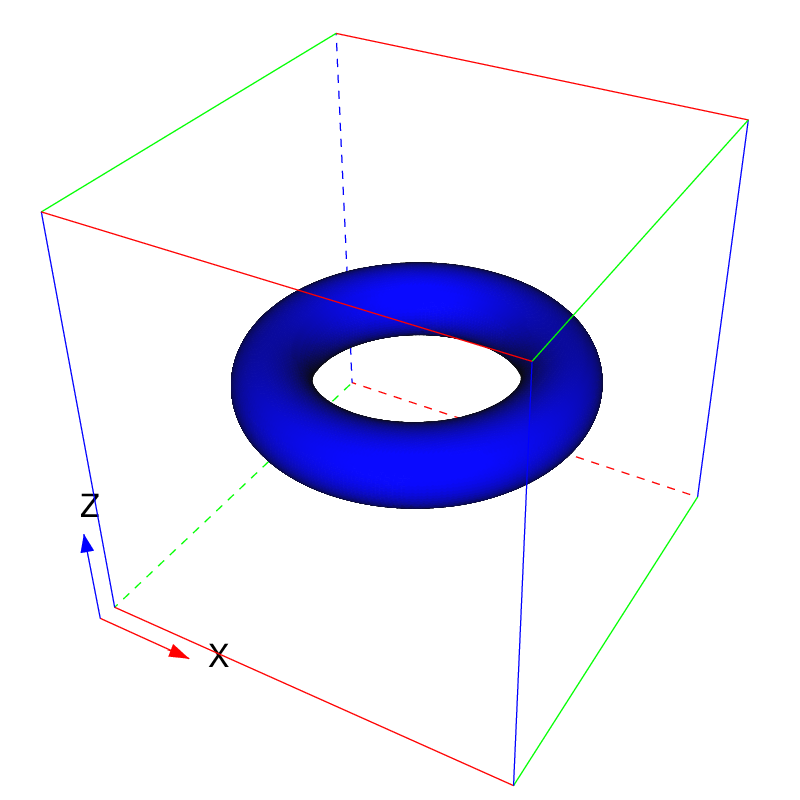}}
\subfigure[t = 0.10 ]{\includegraphics[width=0.3\textwidth]{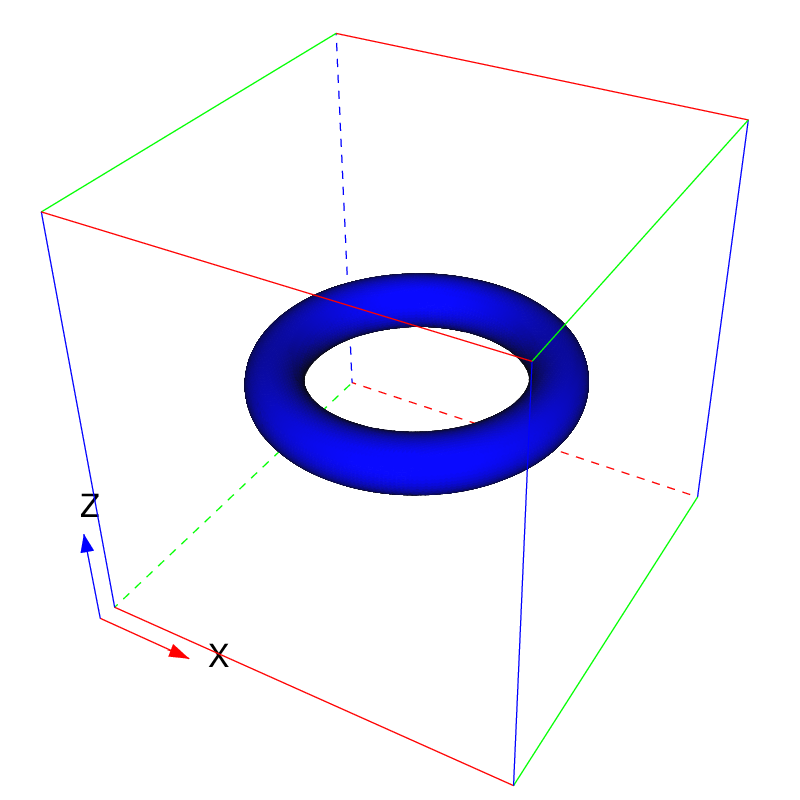}}
\caption{Time evolution of the MCF for a torus-shaped surface.}
\label{fig:torus}
\end{figure}

\begin{figure}[htbp]
    \centering
    \includegraphics[width=0.8\textwidth]{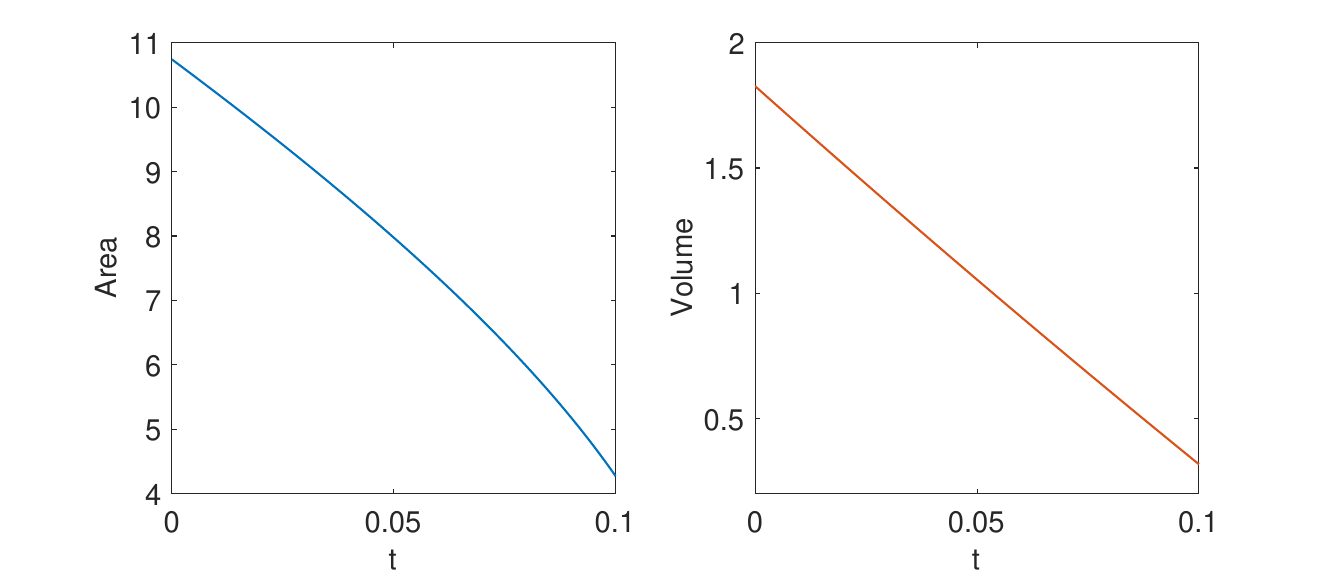}
    \caption{Time evolution of the surface area and enclosed volume of the torus-shaped surface.}
    \label{fig:torus-A-V}
\end{figure}

In the second case, we chose a genus $1$ torus-shaped initial surface.
The surface is given by 
\begin{equation}
\Gamma = \left\{(x,y,z)\Big |\left (c - \sqrt{x^2 + y ^2}\right )^2 + z ^ 2 - a^2=0\right\},
\end{equation}
with $a = 0.34, c = 0.8$.
Numerical error and convergence order are summarized in Table \ref{tab:torus}.
The time evolution of the surface and its area and enclosed volume are presented in Figure \ref{fig:torus} and \ref{fig:torus-A-V}, respectively.
Driven by mean curvature, the torus-shaped surface becomes thinner with time.

\begin{table}[htbp]
\caption{Numerical error and convergence order of the 3D MCF for a torus-shaped initial surface.}
\label{tab:torus}
\centering
\begin{tabular}{|c|c|c|c|c|}
\hline
$N$     & $\Vert \boldsymbol{e}_h\Vert_{\infty}$ & order  &  $\Vert\boldsymbol{e}_h \Vert_{2}$ &  order \\ \hline
128  & 1.65e-03 & - & 3.08e-04 & - \\ \hline
256  & 1.29e-03 & 0.36 & 1.64e-04 &  0.91 \\ \hline
512  & 8.26e-04 & 0.64  & 5.94e-05 &  1.47 \\ \hline
\end{tabular}
\end{table}

\begin{figure}[htbp]
\centering
\subfigure[t = 0 ]{\includegraphics[width=0.3\textwidth]{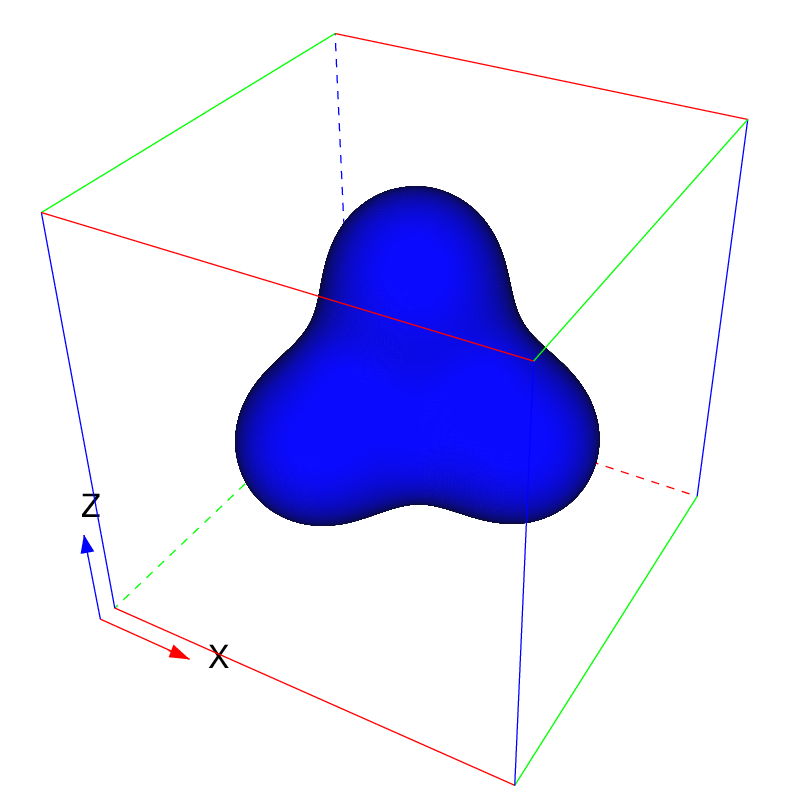}}
\subfigure[t = 0.04 ]{\includegraphics[width=0.3\textwidth]{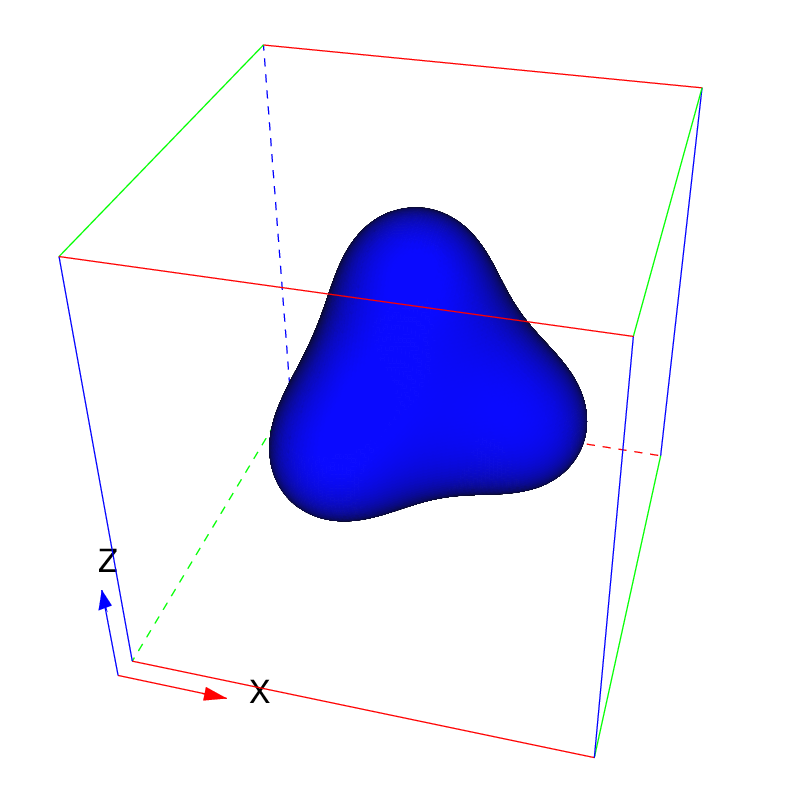}}
\subfigure[t = 0.08 ]{\includegraphics[width=0.3\textwidth]{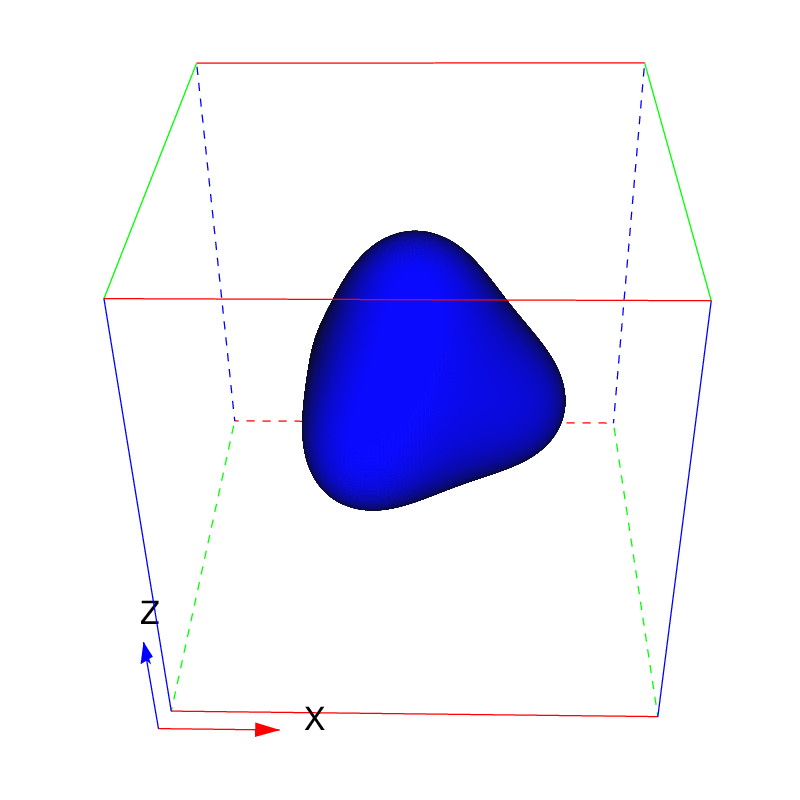}}\\
\subfigure[t = 0.12 ]{\includegraphics[width=0.3\textwidth]{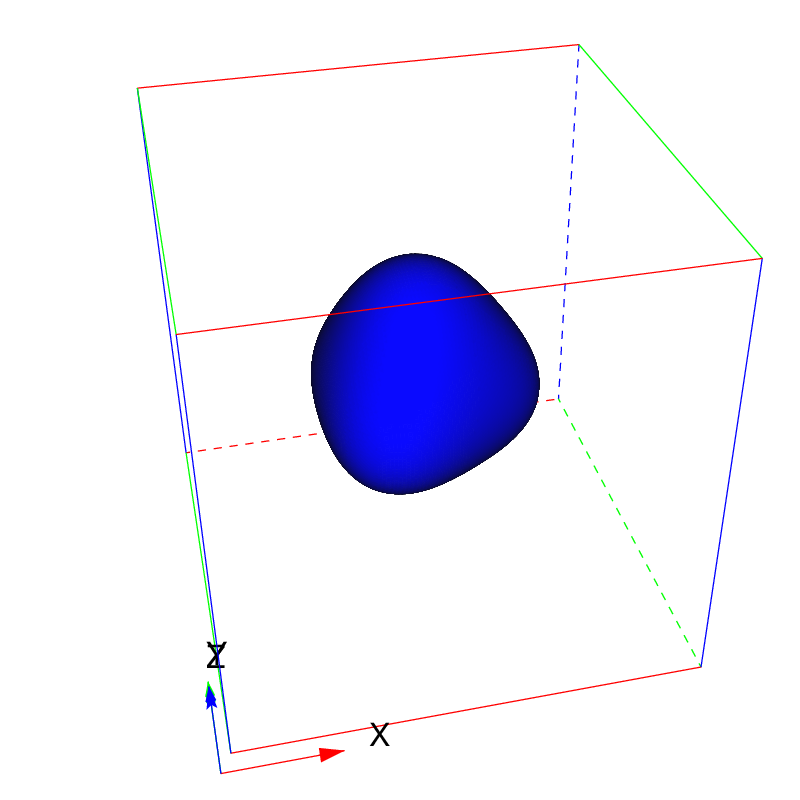}}
\subfigure[t = 0.16 ]{\includegraphics[width=0.3\textwidth]{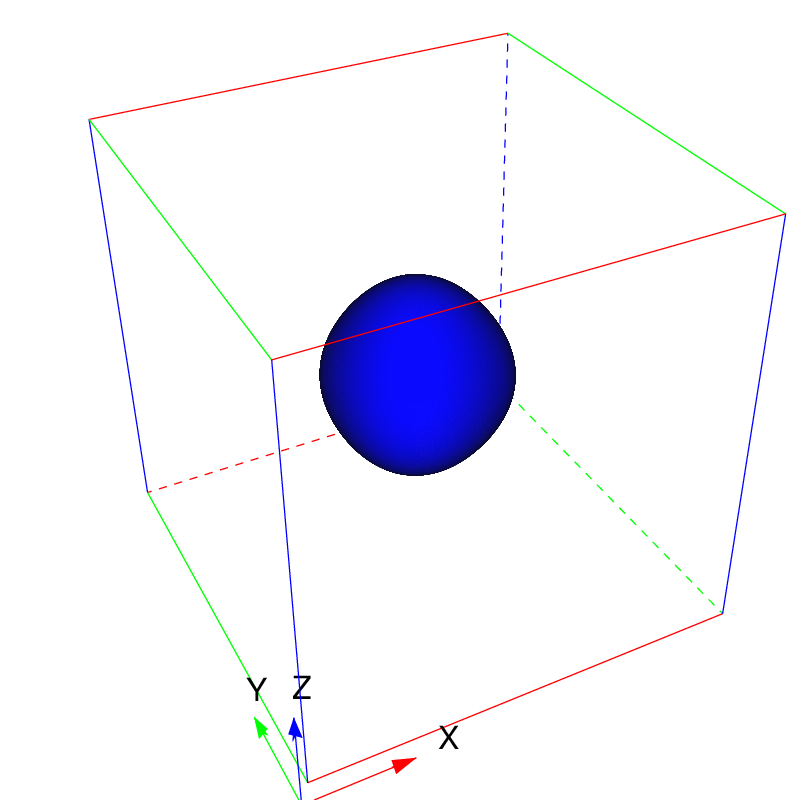}}
\subfigure[t = 0.20]{\includegraphics[width=0.3\textwidth]{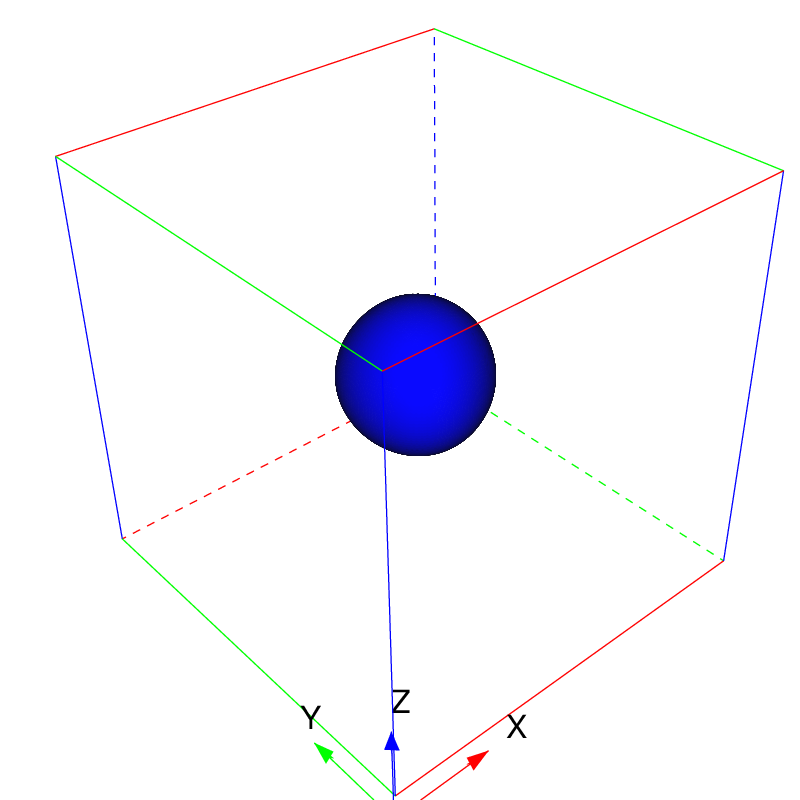}}
\caption{Time evolution of the MCF for a molecular-shaped surface.}
\label{fig:molecular}
\end{figure}

\begin{figure}[htbp]
    \centering
    \includegraphics[width=0.8\textwidth]{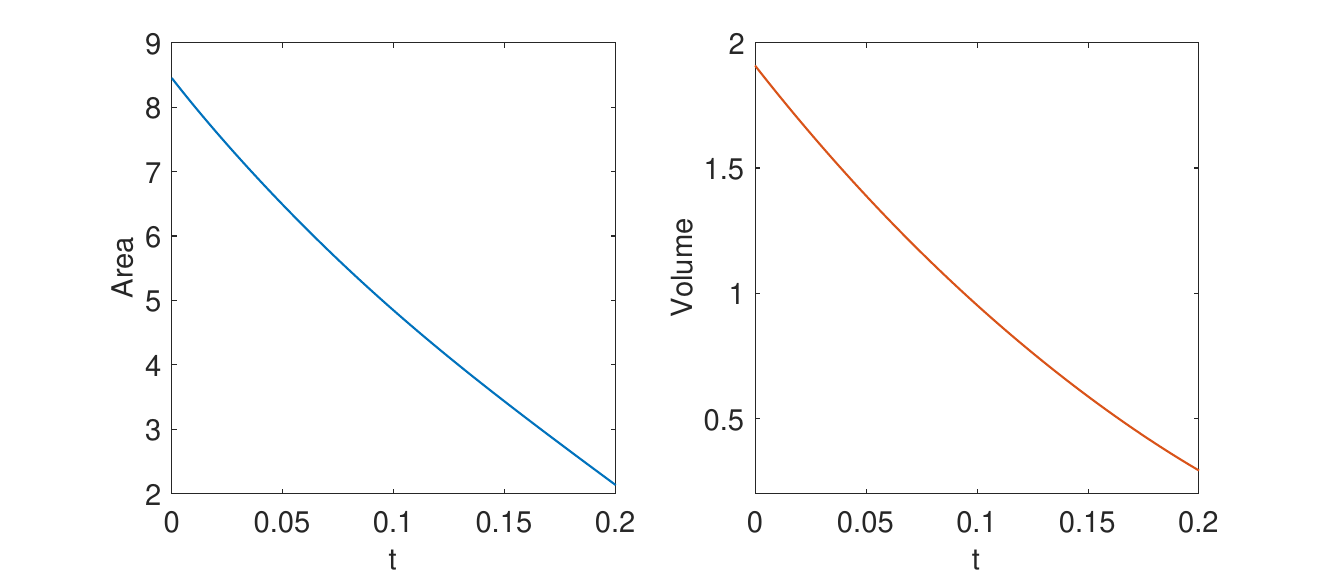}
    \caption{Time evolution of the surface area and enclosed volume of the molecular-shaped surface.}
    \label{fig:atoms-A-V}
\end{figure}

In the final case, the initial surface is a four-atom molecular-shaped surface which is given by
\begin{equation}
\Gamma = \left\{(x,y,z)\Big |c - \sum_{k = 1}^4 \exp\left (-\frac{|\boldsymbol{x} - \boldsymbol{x}_k|^2}{r^2}\right ) = 0\right\},
\end{equation}
with $\boldsymbol{x}_1 = (\sqrt{3} / 3, 0, -\sqrt{6}/12)$, $\boldsymbol{x}_2 = (-\sqrt{3} / 6, 0.5, -\sqrt{6}/12)$, $\boldsymbol{x}_3 = (-\sqrt{3} / 6, -0.5, -\sqrt{6}/12)$, $\boldsymbol{x}_4 = (0, 0, \sqrt{6}/4)$ and $c = 0.5$, $r = 0.5$. 
The numerical error and convergence order are summarized in Table \ref{tab:molecular}.
The time evolution of the surface and its area and enclosed volume are presented in Figure \ref{fig:molecular} and \ref{fig:atoms-A-V}, respectively.

\begin{table}[htbp]
\caption{Numerical error and convergence order of the 3D MCF for a molecular-shaped initial surface.}
\label{tab:molecular}
\centering
\begin{tabular}{|c|c|c|c|c|}
\hline
$N$      & $\Vert \boldsymbol{e}_h\Vert_{\infty}$ & order  &  $\Vert \boldsymbol{e}_h \Vert_{2}$ &  order \\ \hline
128  & 1.79e-03 & - & 4.43e-04 & - \\ \hline
256  & 1.24e-03 & 0.53 & 2.51e-04 &  0.82 \\ \hline
512  & 5.19e-04 & 1.26  & 9.12e-05 &  1.46 \\ \hline
\end{tabular}
\end{table}

\section{Discussion}\label{sec:discu}
This work presents a Cartesian grid-based alternating direction implicit method for solving mean curvature flows in two and three space dimensions.
The method decomposes a hypersurface into multiple overlapping subsets for which new evolution equations are derived by adding extra tangential velocities.
The new formulations for the moving hypersurface only require solving a sequence of scalar quasi-linear parabolic PDEs on planar domains, which is one dimensional lower than the original formulation.
The overlapping subsets of the hypersurface can be represented in terms of height functions of Monge patches which are discretized with Cartesian grids.
With this representation of the hypersurface, an ADI-type semi-implicit time integration method is proposed such that the subsets can be evolved alternately.

The convergence of the proposed method is validated by numerical experiments.
The results show that the ADI method is efficient compared with an explicit scheme since it does not have high-order stability constraints on time step size.
Mean curvature flows for various hypersurfaces in two and three space dimensions are also presented, including one whose initial configuration is a genus $1$ surface.

Although the method in this paper is designed for solving mean curvature flows, it is expected to be able to solve more moving interface problems described by geometric evolution laws, such as the anisotropic mean curvature flow, the surface diffusion flow, and the Willmore flow.
Further, for problems that involve moving interfaces and bulk PDEs simultaneously, such as the Stefan problem and two-phase Stokes flow, the method can also be applicable if combined with a PDE solver such as the kernel-free boundary integral method \cite{Ying2013}.



\textbf{Funding}
W. Y. is financially supported by the National Key R\&D Program of China, Project Number 2020YFA0712000, the Strategic Priority Research Program of Chinese Academy of Sciences (Grant No. XDA25010405), the National Natural Science Foundation of China (Grant No. DMS-11771290) and the Science Challenge Project of China (Grant No. TZ2016002). S. L. is partially supported by the U.S.  National Science Foundation, Division of Mathematical Sciences grants DMS-1720420 and DMS-2309798.

\textbf{Data availibility}
Enquiries about data availability should be directed to the authors.

\section*{Declarations}
\textbf{Conflict of interest} 
We declare that we have no financial and personal relationships with other people or organizations that can inappropriately influence our work. There is no professional or other personal interest of any nature or kind in any product, service and/or company that could be construed as influencing the position presented in, or the review of, the manuscript entitled.

\bibliographystyle{plain}
\bibliography{references}   


\end{document}